\documentclass[draftcls,11pt,onecolumn,a4]{IEEEtran}

\newcommand{\cov}{\ensuremath{\mathrm{Cov}}}
\usepackage{graphics,graphicx,bm,amsmath,amssymb}      
\usepackage{pst-all,amstext,amsfonts,epsf,epsfig,pifont}
\begin{document}
%
\title{Multiple Multidimensional Morse Wavelets\\
{\large Statistics Section\\
Technical Report TR-05-02\\
\today}
}
%
%
\author{Georgios Metikas \& Sofia~C.~Olhede\\
Department of Mathematics, Imperial College London, SW7 2AZ UK
\thanks{Manuscript received XXXXXXX XX, XXXX; revised XXXXXXX XX, XXXX.
G. Metikas was supported by an EPSRC grant.}
\thanks{G. Metikas and S. Olhede are with the Department of Mathematics,
Imperial College London, SW7 2AZ, London, UK (s.olhede@imperial.ac.uk). Tel:
+44 (0) 20 7594 8568, Fax: +44 (0) 20 7594 8517.}}

\markboth{Statistics Section Technical Report TR-05-02}{Metikas
\& Olhede:
Multiple Multidimensional Morse Wavelets}
\maketitle

\begin{abstract}
We define a set of operators that localise a radial image in radial space
and radial
frequency simultaneously. We find the eigenfunctions of this
operator and thus define a non-separable orthogonal set of radial wavelet functions
that may be considered optimally concentrated over a region of radial space and
radial scale space, defined via a doublet of parameters. We give
analytic
forms to their energy concentration over this region. We show how the radial
function
localisation operator can be generalised to an operator, localising any
$L^2({\mathbb{R}}^2)$ function. We show that the latter operator, with an
appropriate choice of localisation region, approximately
has the same eigenfunctions as the radial operator.\\
Based on the radial wavelets we define a set of quaternionic valued wavelet
functions
that can extract local orientation for discontinuous signals and both orientation
and phase structure for oscillatory signals. 
The full set of quaternionic
wavelet functions are component wise orthogonal; hence their statistical
properties are
tractable, and we give forms for the variability of the estimates of the
local
phase and orientation, as well as the local energy of the image. By averaging
estimates across wavelets, a substantial reduction in the variance
is achieved.
\end{abstract}

\section{Introduction}
Localised analyses in one dimension have proven to be remarkably
successful -- notably so wavelet analyses. The latter is based on the
idea that observed signals varying over an increasing argument, time say,
exhibit disparate and highly localised behaviour associated with variations
at a particular scale and at particular time points. Analysis is based on
the wavelet transform, given for signal $g(t)$ using wavelet $\psi(t)$ via
\begin{eqnarray}
\label{wt}
w_{\psi}(a,b;g)&=&\langle \psi_{a,b},g\rangle\\
\nonumber
&=&\int_{-\infty}^{\infty} g(t)\left|a\right|^{-1/2}
\psi^*\left(\frac{t-b}{a}\right)\;dt,
\end{eqnarray}
where $a$ is referred to as the {\em scale}, $b$ the {\em translation} and
$*$ denotes conjugation. This allows for the recognition of patterns specific
to time localisations $b$ and length scales associated with
scale $a,$ if the function $\psi(\cdot)$ is chosen such that the
support
of $\psi(\cdot)$ is essentially limited to a region near the origin, and
the support of the Fourier transform of $\psi(\cdot),$ is essentially limited
to a region near some non-zero
reference frequency $f_{\max}.$ A function cannot be perfectly compact
in time
and frequency simultaneously, and so other criteria have been specified to
determine the localisation properties of $\psi(\cdot).$ Of particular note
is the idea of a localisation operator, generalizing the truncation in time
or frequency
operators \cite{Thomson} to {\em simultaneously} localising in time and frequency/scale
\cite{DaubechiesPaul1988,Olhede2002}. The eigenfunctions/eigenvectors of such operators
are optimally localised with respect to the operator and in one
dimension
the problem of defining appropriate operators and calculating their eigenfunctions
has
been considered in detail \cite{DaubechiesPaul1988,Olhede2002,Lilly}.
\\
The choice of extension of decomposition to
two dimensional analysis is not trivial since variation in the spatial variable
is associated with a direction, as well as a scale. This direction cannot
be assumed to be aligned with the observational coordinate
axes, and thus analysis using a
simple tensor product of one-dimensional wavelets is, in general, not suitable.
In two dimensions
localisation is made to spatial point $\bm{b}=\left[b_1,\;b_2\right]^T,$
in scale to $a$ and in orientation to angle $\theta\in\left[0,2\pi\right)$
{\em cf} \cite{antoinemurenzi99}.
The two dimensional 
continuous wavelet decomposition
of image $g(\bm{x})$ using wavelet $\psi(\bm{x})$ is constructed via
\begin{eqnarray}
\label{wt2d}
w_{\psi}(a,\theta,\bm{b};g)&=&\langle \psi_{a,\theta,\bm{b}},g\rangle\\
\nonumber
&=&\int \int_{\bm{R}^2} g(\bm{x})
\psi_{a,\theta,\bm{b}}^*\left(\bm{x}\right)\;d^2\bm{x},
\end{eqnarray}
where (${\mathcal{D}}_a$) represents a dilation, (${\mathcal{T}}_{\bm{b}}$)
a translation and (${\mathcal{R}}_{\theta}$) a rotation of $\psi(\bm{x})$
giving
\begin{eqnarray}
\psi_{a,\theta,\bm{b}}\left(\bm{x}\right)
&=&{\mathcal{R}}_{\theta}{\mathcal{D}}_a{\mathcal{T}}_{\bm{b}}\psi(\bm{x})\\
&=& \left|a\right|^{-1}
\psi\left(\bm{r}_{-\theta}a^{-1}\left(\bm{x}-\bm{b}\right)\right),
\end{eqnarray}
with $\bm{r}_{\theta}$ given as the rotation matrix.
The decomposition will, with an appropriate choice of wavelet function,
uncover/disentangle behaviour across specific
spatial points, scales and orientation.
Local two dimensional patterns in general
may be intrinsically one dimensional, i.e. after a suitable rotation all
variation is along
a single axis, or intrinsically two dimensional, i.e. there is variation
in several
directions, operating at the same scale emanating from one spatial point.
Following remarks by \cite{Vese} 
we focus on wavelet analysis of discontinuities and oscillatory structure.
Note that edges, or spatial discontinuities,
have an orientation if they locally correspond to (one dimensional) curved
discontinuities
whilst (two dimensional) point discontinuities have no associated orientation.
Oscillations may structurally take the form of one dimensional objects
such
as repeated lines with an even spacing that, if
rotated to the appropriate axes, can be locally described as constant
in one variable and as a sinusoid in the other. Two dimensional
oscillations, circularly emanating
from a single point, when considered locally at a distance from their
source
may be described approximately as one dimensional oscillations.
\\
Thus the structure of an image is highly orientation dependent, and
analysis methods should disentangle both locally one dimensional and two
dimensional
structures operating at many different orientations. 
A well known feature of wavelet analysis, \cite{Donoho}
is that genuinely two dimensional structure, {\em i.e.}
point discontinuities, are well represented in a wavelet decomposition, and
so in this paper
we focus instead on the treatment of locally one dimensional structures by
adjusting the wavelet transform suitably. We shall discuss
two existing strategies for considering oriented scale-based decompositions,
and 
construct a new method corresponding to
a synthesis of the two methods discussed. This method
can extract the local orientation of the image explicitly, in a multi-scale
framework.\\
Existing continuous wavelet methods that deal with the orientation of the
image explicitly
are based on directionally selective filters, or 
{\em directional wavelets.}
Antoine \& Murenzi \cite{antoinemurenzi99} define complex
directional wavelets with
a preferred orientation in the frequency domain, as their frequency
support is limited to a pre-defined cone, parameterised
via the opening and closing angles of the cone \cite[p.~324--5]{antoinemurenzi99}. 
Defining highly directional
wavelets will necessitate an elongation of the wavelets in the spatial frequency
domain, and this affects their
spatial and spatial frequency resolution capacity -- along a specific orientation
in the frequency domain, the wavelets localise badly in frequency, and for
this reason we do not use directional wavelets.
Directional wavelets will localise in scale {\em and} direction, whilst spatially
isotropic wavelets only separate features at different scales. 
We shall 
use an isotropic wavelet decomposition to separate
out disparate components occurring at either different scales and the same spatial
locations, or at the same scales at different spatial locations.
To facilitate this separation of structure, wavelets that are optimally concentrated
in radial space and radial frequency are required, as in two dimensions the
notion
of spatial
distance is naturally associated with the Cartesian metric. We
define a family of radial two-dimensional localisation operators and find
the
radial eigenfunctions
of any given operator in this family, denoted the isotropic Morse wavelets. Any operator in this family is characterised
via 
two parameters that determine the spatial/spatial frequency structure of
the isotropic eigenfunctions.
Any choice of the parameters fixes a particular operator that in turn possesses
a family of eigenfunctions. These functions
are orthogonal and indexed via $n.$ The eigenvalues explicitly give the radial
concentration of the eigenfunctions. These eigenfunctions are related (but
not equivalent) to the eigenfunctions of the one dimensional
Morse localisation operator \cite{Olhede2002}.
For every fixed value of $n,$ and given radial eigenfunction, we define an
additional pair of functions,
whose joint norm may be considered to have the same localisation in space
and spatial scale as the original radial function, but when combined with
the original radial function will disentangle the local orientation
of the image
analysed. These extra pairs of functions are constructed explicitly to consider
local orientation and phase.\\
A method for considering local phase
structure is to extend the notion of instantaneous frequency 
\cite{Boashash1992}; \cite{FelsbergSommer,vonBulow,Hahn1992} give extensions
to
instantaneous frequency and local phase structure in the spatial domain,
and for each spatial point
retrieve a local phase/variational structure. These extensions
correspond to the calculation of several additional images, or quadrature
components, at each spatial point, 
where each additional set of components is considered to have the same local spatial
energy and variational structure as the original image. The full set of components is used to
calculate the local orientation and variational structure.
The additional components are ill-defined when constructed from multi-component
images, as then
a single component with a spatially varying phase function is not an appropriate model for
the original image. It thus becomes necessary to combine the calculation of a local
phase with scale-localised
methods such as the wavelet transform.
For each fixed value $n\in{\mathbb{N}}$, for each Morse
wavelet, we define two extra real functions
to complement the isotropic Morse wavelets. The triplet of real valued functions
form a monogenic signal \cite{FelsbergSommer}. Two of the triplet of functions
should
be thought of as a single vector valued object, where their vector structure characterises the orientation of the local variations. 
The monogenic wavelets and wavelet transform are best represented
using quaternion \cite{Deavours}, rather than, real or complex numbers. Each
triplet is therefore considered as a positive real valued amplitude, a
pure unit quaternion specifying a direction, and a phase. The real amplitude characterises
a local energy, the pure unit quaternion an orientation, and the phase a
local variational structure \cite{FelsbergSommer}.
The quaternion algebra allows for easy parameterisation of phase and orientation
structure. 
\\
In one dimension analytic wavelets, {\em i.e.}
Cauchy or
Morlet wavelets \cite[p.~28]{Holschneider1998}, are used
to identify local oscillatory structure, of a real image $g(\cdot).$ 
The monogenic Morse wavelets are the natural two-dimensional
extension of the analytic Morse wavelets \cite{Olhede2002}, and they define
a local phase
and orientation structure at each spatial/spatial scale point for a
two dimensional oscillatory image, yielding a natural, and more elegant, structure for
wavelet ridge analysis \cite{Gonnet}. Oscillatory images are the complements
of discontinuous images, where
the two kinds of images appear as either oscillatory or discontinuous depending on if they
are analysed in the Fourier or spatial domain.
We consider line, or curved, discontinuities that may locally be thought
of
as one dimensional structures, and discuss the retrieval of their local features.
Analyses using the separable discrete wavelet transform in addition with local phase structure characterisations, have been previously considered \cite{Chan,Fernandes}.
However, in contrast to their procedure,
we define {\em multiple orthogonal continuous} wavelets, based on
a different
two dimensional extension to the analytic signal. The Morse wavelets are
additionally
optimally concentrated with respect to a radial position/scale region ${\mathcal{D}}.$ Other quaternion valued decompositions
includes the work of \cite{Hahn2005}, however this decomposition is only
suitable for deterministic images.
\\
We briefly discuss the discrete implementation of the two dimensional monogenic
Morse
wavelet transform, and in more depth the statistical properties of the transform.
Unavoidably, most observed signals are contaminated by noise, and so robust
methods that can deal with noise must be designed. As the operator problem
yields solutions of multiple
orthogonal monogenic wavelets, we may use the notion of averaging
uncorrelated estimates \cite{Thomson} to retrieve estimates of the space-scale
energy as well as other quantities of the image with reduced variance. Multiple
orthogonal
filters have been considered in several dimensions for stationary processes
\cite{deville}, and non-stationary processes \cite{Simons2003} using the
windowed Fourier
transform and tensor product windows, but our wavelets 
are in contrast to these methods orthogonal, {\em non-separable monogenic}
wavelet functions.
Finally, the methods are illustrated on typical examples, showing the power
of
the multiple monogenic Morse wavelets.

\begin{figure*}[t]
\centerline{
\includegraphics[height=1.70in,width=1.70in]{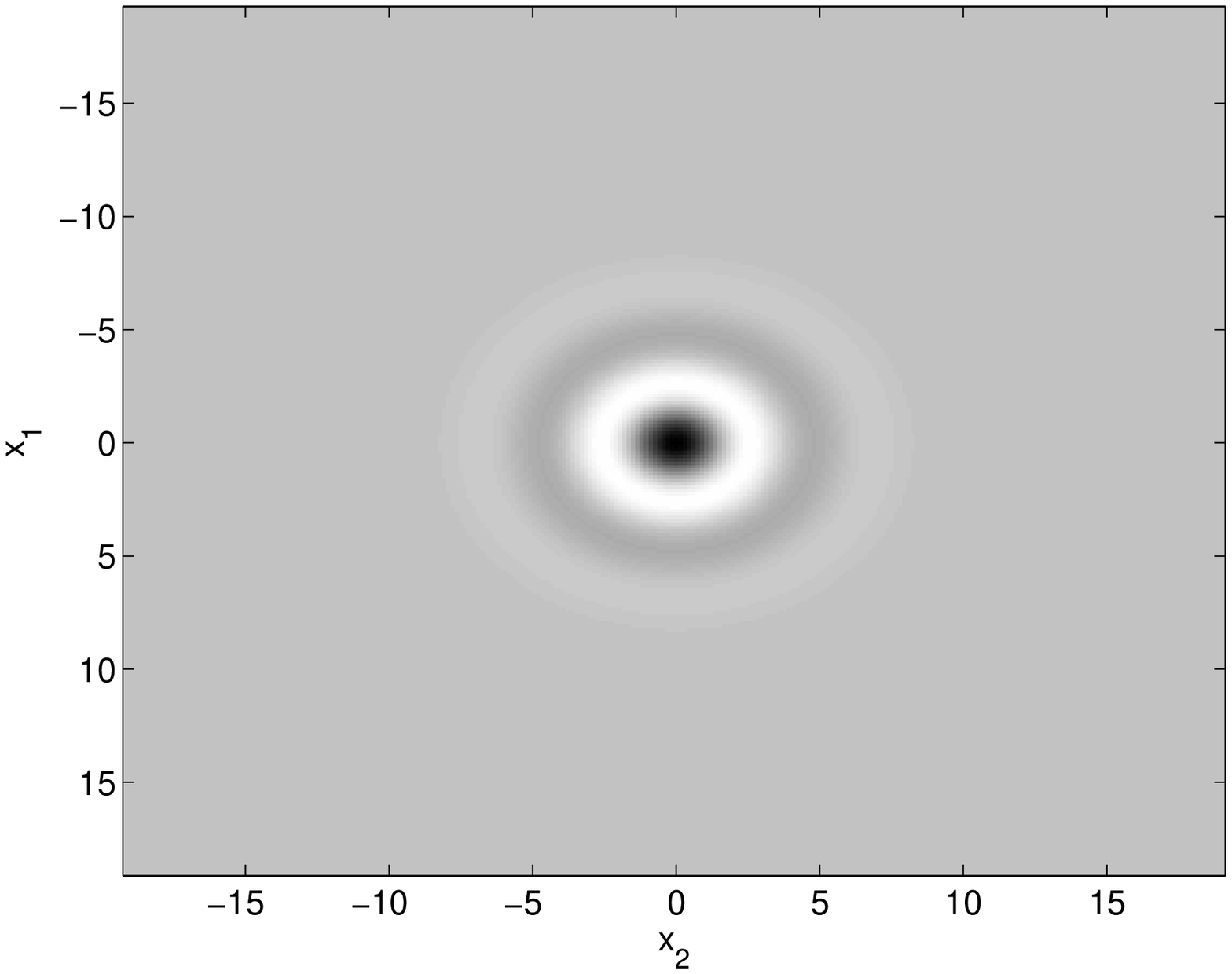}
\includegraphics[height=1.70in,width=1.70in]{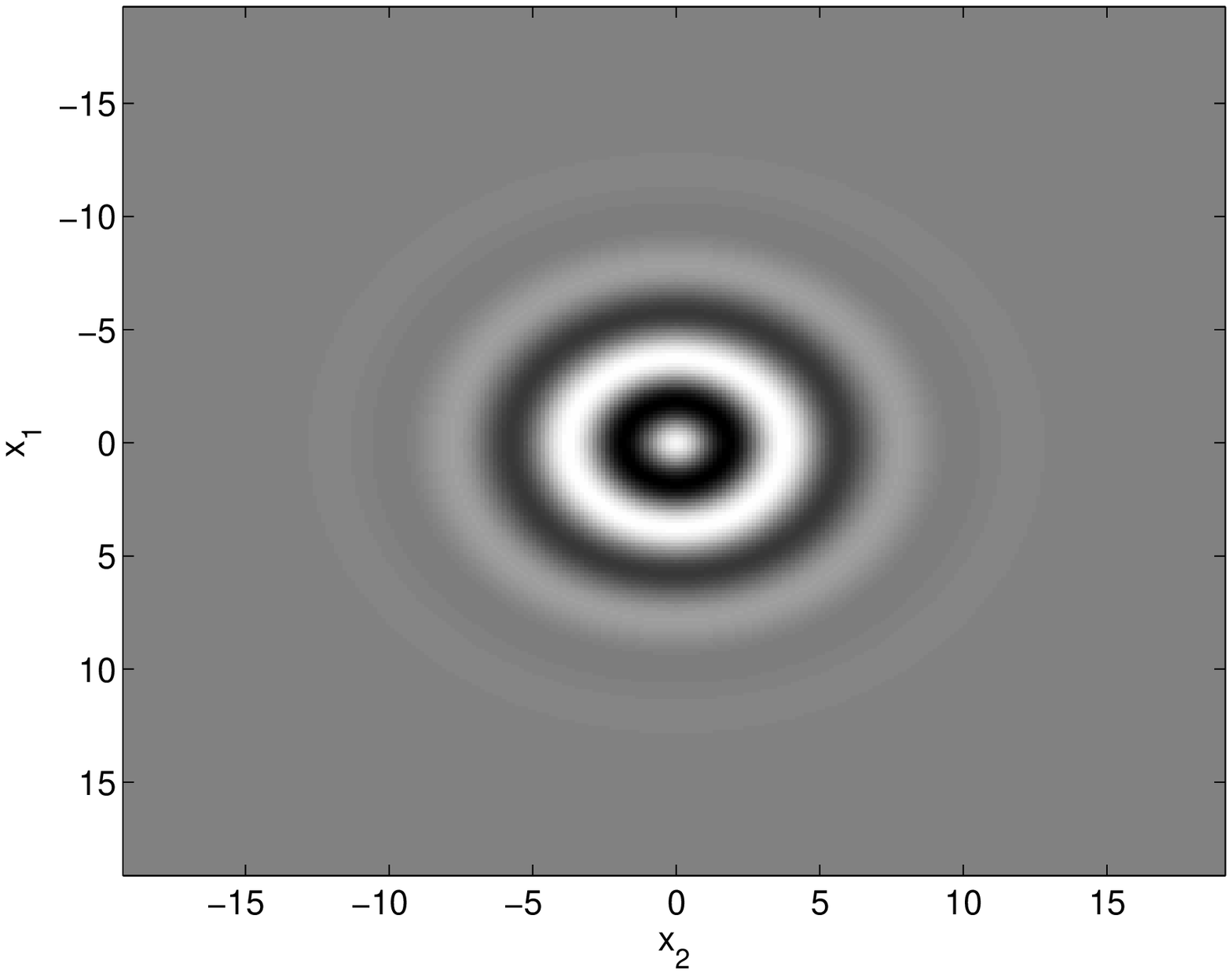}
\includegraphics[height=1.70in,width=1.70in]{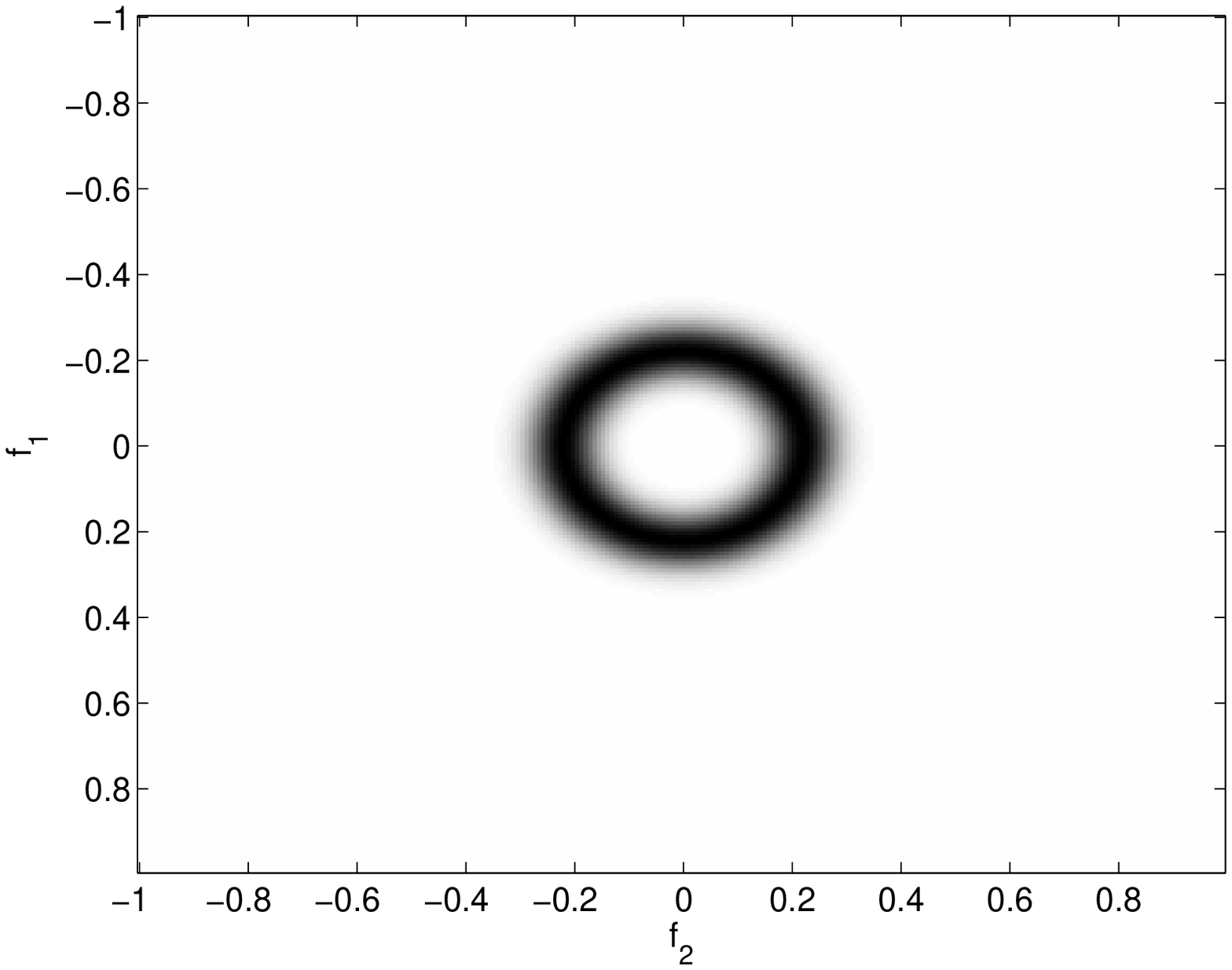}
\includegraphics[height=1.70in,width=1.70in]{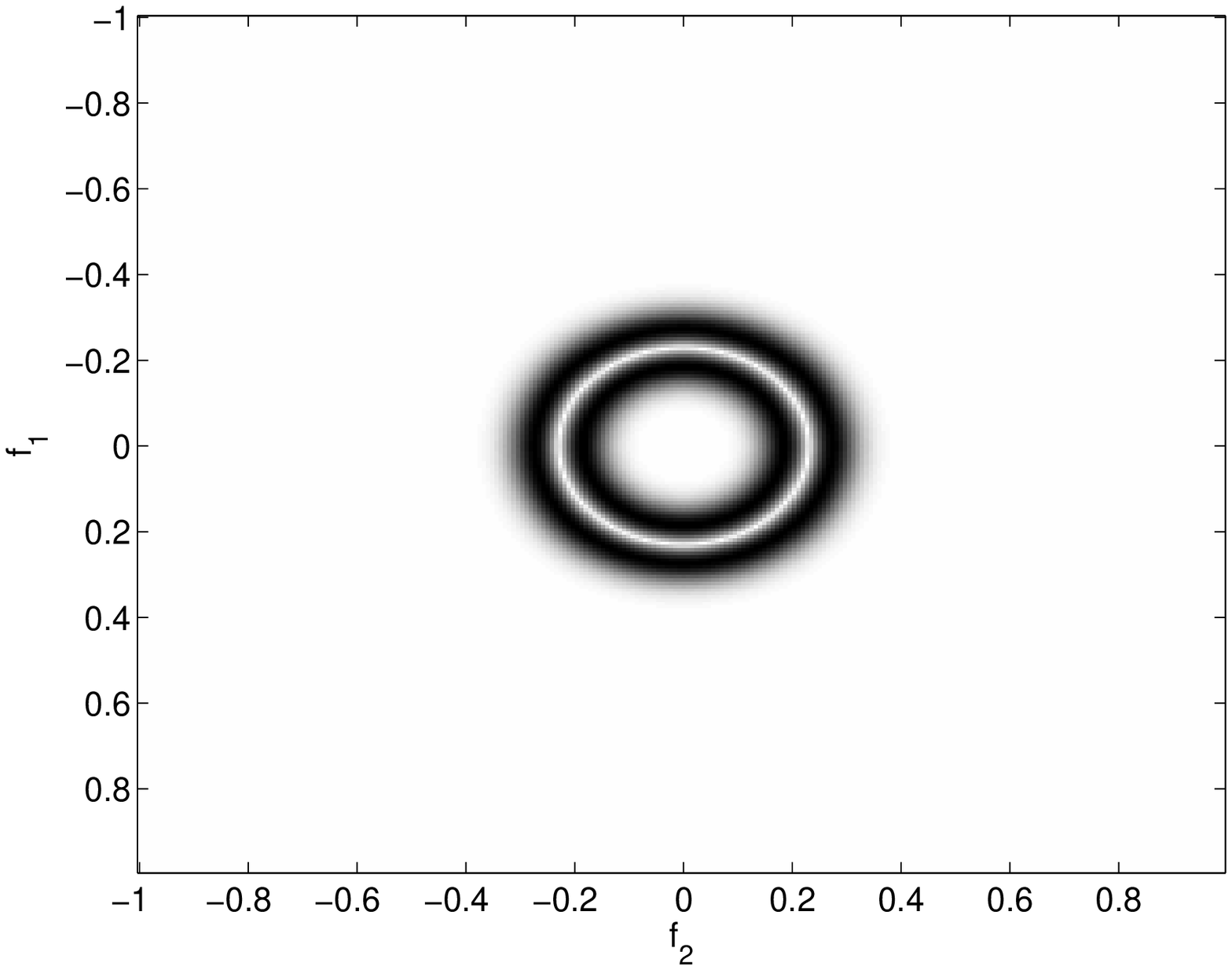}}
\caption{\label{fig:partition} 
The isotropic Morse wavelets in the spatial domain.
$(l,m)=(8,3)$ and $n=0$ (far left), $n=1$ (second from left). The modulus
of the isotropic Morse wavelets in the spatial frequency domain.
$(l,m)=(8,3)$ and $n=0$ (second from right), $n=1$ (far right).
}
\end{figure*}

\section{Notation}
We denote the 1D Fourier transform of $g(\cdot),$ $G(f)=\langle
\exp(2\bm{j}\pi ft),g\rangle,$ and the two dimensional Fourier transform of $g(\bm{x})$
as $G(\bm{f})=\langle
\exp(2\bm{j}\pi \bm{fx}),g\rangle.$ 
We denote an arbitrary quaternion via $q=q_1+q_2
\bm{i}+q_3\bm{j}+q_4 \bm{k},$ where $q_i\in{\mathbb{R}},\;i=1,\dots4$ and
note that $\bm{i}^2=\bm{j}^2=\bm{k}^2=\bm{i}\bm{j}\bm{k}=-1,$
whilst $\bm{ij}=-\bm{ji}=\bm{k},$ $\bm{ik}=-\bm{ki}=-\bm{j},$ 
and finally $\bm{jk}=-\bm{kj}=\bm{i}.$
The algebra is non-Abelian.
We additionally define the two dimensional Fourier transform in terms of any unit quaternion as
$G_{\bm{e}_q}(\bm{f})=\langle
\exp(2\bm{e}_q\pi \bm{fx}),g\rangle,$ so that the regular Fourier transform
corresponds to $G_{\bm{j}}(\bm{f})\equiv G(\bm{f}).$ 
For more
notes on quaternion algebra see \cite{Deavours}. We retain here only the
briefest
possible usage of the quaternion algebra, necessary for clarity of exposition,
and stress that all implementation
is discussed in terms of real vector quantities.
The rotation
operation is implemented using matrix 
\[r_{\theta}\equiv
\begin{pmatrix} \cos(\theta) & -\sin(\theta)\\
\sin(\theta) & \cos(\theta)
\end{pmatrix}.\]
We shall extensively
use polar coordinates, and define: $\bm{x}=\left[x \cos(\chi),
x \sin (\chi)\right],$ $\bm{f}=\left[f \cos(\phi),f \sin (\phi)\right],$
and $\bm{b}=\left[b \cos(\phi_b),
b \sin (\phi_b)\right].$

\section{Two dimensional Wavelet Analysis}
One dimensional local analysis corresponds to decomposing a function in terms
of
a set of functions that contain behaviour local to a set scale $a,$ and time
point $b.$
The two dimensional wavelet transform is defined using four parameters:
$\bm{\xi}
=\left[a,\theta,\bm{b}\right],$
where $a$ and $\bm{b}$ play roughly the same role as the corresponding
one dimensional quantities, and $\theta$ corresponds to local orientation
localisation.
For the decomposition to be meaningful, the translated and dilated wavelets
are chosen in one dimension to be mainly supported near time point $b$ and frequency
point $f_a=f_{\max}/a,$ and the obvious extension to two dimensions would be to
find functions that are mainly supported at spatial point $\bm{b}$ and spatial
frequency point $\bm{f}_a=\frac{f_{\max}}{a}\left[\cos(\theta),\;\sin(\theta)\right].$
To measure the localisation of a given wavelet function, commonly
its spread in time and frequency or scale is calculated, and both quantities
are desired to be low -- their product is bounded below, thus restricting
possible joint localisation.
A function's spread in time and frequency should be considered
simultaneously in the two domains
\cite{Olhede2002},
rather than combining two separate marginal properties: for this purpose
localisation operators were defined, denoted ${\mathcal{P}}_{\mathcal{D}}.$
The localisation of a function $g(\cdot)$ to region ${\mathcal{D}}$ is measured
by the ratio of energy $\mu_g({\mathcal{D}})=\langle {\mathcal{P}}_{\mathcal{D}}g(\cdot),
{\mathcal{P}}_{\mathcal{D}}g(\cdot)\rangle/\langle g(\cdot),g(\cdot)\rangle,$
and the eigenfunctions of ${\mathcal{P}}_{\mathcal{D}}$ achieve optimal ratios
\cite{Olhede2002}.
Naturally extending analysis to two dimensions requires the appropriate definition
of a two-dimensional localisation operator. We consider localisation
in scale and spatial location, which requires localisation in radial spatial
frequencies of spatially radial
functions, as we wish to associate a notion of distance to the Cartesian
two dimensional metric. This produces functions optimally concentrated
in a radial space and scale $a,$ but that have no orientation. We
then construct a set of functions capable of extracting local orientation
information
based on the monogenic signal \cite{FelsbergSommer}. \cite{Dahlke} have shown
that with an appropriate definition
of scale, orientation and spatial spread, the orientation and scale may
be considered separately when finding optimally concentrated functions.
In a slightly different setting they demonstrate that radial functions have
optimal scale versus position
localisation
properties, for appropriately chosen families of functions. This, in combination
with a Cartesian metric in the plane, motivates the study of radial localisation.

\subsection{Radial Localisation Operator}
The construction of a coherent radial state, required for the construction
of a radial
projection operator, is not straightforward. 
We commence with radial function $v(\bm{x})=v_r(x)$ and in theory wish to
construct a family of radial functions that have been shifted in scale
and position. Obviously this is impossible, as once we shift in position,
we no longer have a radial function, but we may relax our requirements as
the family only needs to act in a similar fashion to a direction averaged
spatial
shift, based on an appropriate domain ${\mathcal{D}}$ chosen. To construct
a generic two-dimensional non-radial coherent
state we would use the full set of parameters $\bm{\xi}=\left[a,\theta,\bm{b}
\right],$ but intend to use the sub-set $\bm{\xi}_{r}=\left[a,b\right].$ Working
with radial functions the act of rotation, represented
by $\theta$ will have little importance, and is not included. The dilation
by $a$ will act correspondingly in the world of radial operations, to one
dimensional analysis, and a substitute for a shift in position by $\bm{b}$
must be defined.
We reconstruct the function using building blocks of $\langle v_{\bm{\xi}},g\rangle
v_{\bm{\xi}}$ that are averaged over suitably defined domains. The localisation domain will be defined in terms
of $b$ and range over $0\le \phi_{b}\le 2\pi.$ This implies that any
shift in $\bm{b}$ will be averaged over the full range of $\phi_{b}$ and
so rather than multiplying $V(\bm{f})$ by $e^{-2\pi \bm{j} \bm{f}\cdot \bm{b}}$
we could multiply this by its angular average of $J_0(2\pi f b),$ that after
implementing the operator projection, yields the same results as the former
strategy. Unfortunately,
this choice of coherent state leads to a mathematically intractable operator,
and
so we define a radial coherent state for $a>0$ and $b>0$ that has similar
properties to the suggested
state via
\begin{equation}
V^{(+)}_{\bm{\xi}_r}(\bm{\omega})=\sqrt{2} a^{1/\gamma}(
a^{1/\gamma} \omega)^{-1/2+\gamma/2} V^{\beta,\gamma}
(a^{1/\gamma}\omega)\frac{\cos\left(\omega^{\gamma}b-\pi/4\right)}
{\sqrt{\omega^{\gamma}
b}}
\label{cohhy2},
\end{equation}
based on the Morse coherent state \cite{Olhede2002} of 
\[
V^{\beta,\gamma}(\omega)
=\frac{2^{r/2+1/2}\sqrt{\pi\gamma}}{\sqrt{\Gamma(r)}}
\omega^{\beta} e^{-\omega^{\gamma}}\;\;{\mathrm{if}}\;\;\omega>0.
\]
Note that $V^{\beta,\gamma}(\omega)\propto \omega^{\beta},\;\beta>0$ for $\omega<<1,$
and for future reference define $\omega_{\epsilon}:\;V^{\beta,\gamma}(\omega)<\epsilon\;\forall
\omega<\omega_\epsilon.$
In equation (\ref{cohhy2}) the dilation of $a$ is implemented as
in \cite{Olhede2002} and requires no further discussion. The factor $
(
a^{1/\gamma} \omega)^{-1/2+\gamma/2}$ is added to ensure the correct normalisation
of the two-dimensional radial function, as is the replacement of $a^{1/\gamma}$
for $a^{1/(2\gamma)},$
as the term in $a$ multiplying the dilated function of $\omega.$ 
Denote the translation-like operator, and the originally posited direction
averaged
operator respectively as,
$
T_{b}^{(1)}(\omega,b)=\frac{\cos\left(\omega^{\gamma}b-\pi/4\right)}
{\sqrt{\omega^{\gamma}
b}},$ and
$T_{b}^{(2)}(\omega,b)=
J_0(\omega^{\gamma}b).$
The translation
operator's decay for large
values of $\omega^{\gamma}
b$ is the same as $T_{b}^{(2)}(\omega,b),$ as is its zero crossing structure,
and apart from small arguments the two functions are performing a similar
action.
The functional behaviour for small values of $\omega$ is different
--  $T_{b}^{(2)}(\omega,b)=1+O\left(\omega^{2\gamma} \right)$ whilst 
$T_{b}^{(1)}(\omega,b)=\frac{1}{\sqrt{2\omega^{\gamma} b}}+O\left(\omega^{\gamma/2}\right),$
where the latter is unbounded near $\omega \rightarrow 0.$
If we denote the coherent state as
\[
V^{+(u)}_{\bm{\xi}_r}(\bm{\omega})=
\sqrt{2} a^{1/\gamma}(
a^{1/\gamma} \omega)^{-1/2+\gamma/2} V^{\beta,\gamma}
(a^{1/\gamma}\omega)
T_{b}^{(u)}(\omega,b),\;u=1,2, \]
then for small values of $\omega,$ with $r=(2\beta+1)/\gamma,$ we find
\begin{eqnarray}
\nonumber
V^{+(u)}_{\bm{\xi}_r}(\bm{\omega})&=&
\sqrt{\frac{2^{r+2}\pi \gamma}{\Gamma(r)}} a^{1/\gamma}(
a^{1/\gamma} \omega)^{-1/2+\gamma/2} 
(a^{1/\gamma}\omega)^{\beta}\left[1+O\left(\omega^{\gamma}\right)\right]
T_{b}^{(u)}(\omega,b)\\
\nonumber
V^{+(1)}_{\bm{\xi}_r}(\bm{\omega})
&=&C_1\omega^{\beta-1/2}+O\left( \omega^{-1/2+\gamma+\beta}\right)\\
V^{+(2)}_{\bm{\xi}_r}(\bm{\omega})
&=&C_2
\omega^{-1/2+\gamma/2+\beta}+O\left(\omega^{\gamma}\right),
\nonumber
\end{eqnarray}
so that both choices give contributions of negligible magnitude for 
$\left|\omega\right|\rightarrow
0$ as long as we assume $\beta>1/2,$ which is combined with the previous
constraint of $\gamma \ge 1,$
with $\beta>(\gamma-1)/2$ \cite[p.~2663]{Olhede2002}.
We choose $T^{(1)}_b$ as the radial spatial
shift quantity. This has an approximate interpretation of $b$
as a orientation averaged spatial shift. Note that the generalised Morse
wavelets also are
based on the approximate, rather than exact, notion of a {\em warped}
location shift \cite[p.~2663]{Olhede2002}. 

We define the operator for a radial function $g(\bm{x})$ in terms of radial
inner product $\langle g_1 ,g_2\rangle_R=\frac{1}{2\pi}
\langle G_1 ,G_2\rangle_R=\int G_1^*(f)G_2(f) f\;df$
as and (in terms of
${\mathcal{D}}^+(C)=\left\{(a,b):\;a^2+b^2+1\le 2 aC,b>0\right\}$),
\begin{eqnarray}
{\mathcal{P}}_{\mathcal{D}^+}\left\{ G \right\}(\frac{\bm{\omega}_1}{2\pi})&=& C_{o}  \int_{ \mathcal{D}  } V^+_{ \bm{\xi} }(\bm{\omega_1}) \langle V^+_{ \bm{\xi} }, G  \rangle  \frac{da}{a^2} db \nonumber \\ 
&=&C_o \int_{\omega_2=0}^{\infty} \int \int_{a^2+b^2+1\le 2a C}\kappa_2(\omega_1,\omega_2;a,b)
\;\frac{da}{a^2}\;db\;d\omega_2
\label{oppy},
\end{eqnarray}
where 
\begin{equation}
\kappa_2(\omega_1,\omega_2;a,b)=2 a^{1/\gamma}
V^{\beta,\gamma}(a^{1/\gamma}\omega_1) \sqrt{\omega_2/\omega_1}
V^{\beta,\gamma}(a^{1/\gamma}\omega_2) G(\frac{\omega_2}{2\pi}) \cos((\omega_1^{\gamma}- \omega_2^{\gamma})b).
\label{kappa2}
\end{equation}
The kernel $\kappa_{2}(\cdot,\cdot;\cdot,\cdot)$ can 
be inverted in its first argument to retrieve the spatial domain operator.
Note that by definition $\omega_1>0.$
We constrain the norm of the operator -- exactly reconstructing the entire radial function if we let the region of integration
across $\bm{\xi}_r$
be large enough, thus calibrating the operator to make the eigenvalues
meaningful.
\begin{eqnarray}
\nonumber
G(\frac{\omega_1}{2\pi})&=&C_o
\int_{\omega_2=0}^{\infty} \int_{a=0}^{\infty} \int_{b=0}^{\infty}
2 a^{1/\gamma}
V^{\beta,\gamma}(a^{1/\gamma}\omega_1) \sqrt{\omega_2/\omega_1}
V^{\beta,\gamma}(a^{1/\gamma}\omega_2) \\
\nonumber
&&G(\frac{\omega_2}{2\pi}) \cos((\omega_1^{\gamma}- \omega_2^{\gamma})b)
\;\frac{da}{a^2}\;db\;d\omega_2\\
\nonumber
&=& C_o \int_{\omega_2=0}^{\infty} \int_{a=0}^{\infty} \int_{b=-\infty}^{\infty}
 a^{1/\gamma}
V^{\beta,\gamma}(a^{1/\gamma}\omega_1) \sqrt{\omega_2/\omega_1}
V^{\beta,\gamma}(a^{1/\gamma}\omega_2) \\
\nonumber
&&G(\frac{\omega_2}{2\pi}) \cos((\omega_1^{\gamma}- \omega_2^{\gamma})b)
\;\frac{da}{a^2}\;db\;d\omega_2\\
\label{cnought}
&=&C_o 2 (2\pi)^2 G(\frac{\omega_1}{2\pi})\frac{1}{r-1},
\end{eqnarray}
which gives $C_o=\frac{r-1}{2^3 \pi^2},$ yielding a  `resolution of identity' \cite{Olhede2002} for radial functions: this will {\em not} hold for
any $g(\cdot) \in L^2(\mathbb{R}^2),$  as the operator is {\em only} defined
for radial functions.
The eigenfunctions of the operator defined in (\ref{oppy}) can be found by solving the equation
\begin{equation}
\label{eigen1}
{\mathcal{P}}_{\mathcal{D}^+}\left\{G\right\}(\frac{\bm{\omega}_1}{2\pi})=\lambda G (\frac{\bm{\omega}_1}{2\pi}).
\end{equation}
The one-dimensional Morse wavelet projection operator can be considered
to be in the case of $\psi(t)$ real, assuming that $\omega_1>0,$
(similar expressions are derived for $\omega_1<0,$ but as we shall use this
to obtain solutions to equation (\ref{eigen1}) we only need to consider $\omega_1>0,$
and the $\sin(\cdot)$ term vanishes due to the symmetry of the projection
region, and where ${\mathcal{D}}$ is as in \cite{Olhede2002}):
\begin{eqnarray}
{\mathcal{P}}_{\mathcal{D}}^{1D}\left\{\Psi\right\}(\frac{\omega_1}{2\pi})&=& C_o^{1D} \int_{\omega_2=0}^{\infty } \int \int_{a^2+b^2+1\le 2a C,\;a\in R^+,\;b\in R} a^{1/\gamma}  V^{\beta,\gamma}
(a^{1/\gamma}\omega_1)V^{\beta,\gamma*}(a^{1/\gamma}\omega_2) \\ \nonumber
&& \Psi(\frac{\omega_2}{2\pi}) \cos{\left[ b (\omega_1^{\gamma}-\omega_2^{\gamma})\right]}
\;\frac{da}{a^2}\;db\;\frac{d\omega_2}{2\pi}\\
\nonumber
&=& \frac{ C^{1D}_{o }}{2\pi} \int_{\omega_2 =0 }^{\infty} \int \int_{a^2+b^2+1 \leq 2aC, \; a\in R^{+}, \; b \in R^{+} } \kappa_{1}(\omega_1,\omega_2,a,b)
 \;\frac{da}{a^2}\;db\;d\omega_2,
\label{ref2}
\end{eqnarray}
where 
\begin{equation}
\kappa_{1}(\omega_1,\omega_2;a,b) = 2 a^{1/\gamma}
V^{\beta,\gamma}(a^{1/\gamma}\omega_1) 
V^{\beta,\gamma}(a^{1/\gamma}\omega_2) \Psi(\frac{\omega_2}{2\pi}) \cos((\omega_1^{\gamma}-
\omega_2^{\gamma})b).
\label{kappa1}
\end{equation}
Note that the above kernel $\kappa_{1}(\cdot,\cdot;\cdot,\cdot)$ is similar
to the kernel $\kappa_{2}(\cdot,\cdot;\cdot,\cdot)$ of (\ref{kappa2}), the
only difference being that $\kappa_2(\cdot,\cdot;\cdot,\cdot)$ has the extra
term, $\sqrt{\omega_2/\omega_1},$ and $C_o^{1D}=\frac{r-1}{4\pi}.$
The Morse wavelets \cite{Olhede2002} are the solution to the equation
\begin{equation}
\label{eigeny2}
{\mathcal{P}}_{\mathcal{D}}^{1D}\left\{\Psi\right\}(\frac{\omega_1}{2\pi})=\lambda^{1D}
\Psi(\frac{\omega_1}{2\pi}).
\end{equation}
Consider equation (\ref{eigen1}); multiply both sides by $\sqrt{\omega_{1}/(2\pi)}$, set $G(\omega)\sqrt{\omega}=\Psi(\omega)$, and note that 
the equation to be solved has now exactly the form of
equation (\ref{eigeny2}):
\begin{eqnarray}
\sqrt{\frac{\omega_1}{2\pi}} {\mathcal{P}}_{\mathcal{D}^+}
\left\{ G \right\}(\frac{\bm{\omega}_1}{2\pi})
=  {\mathcal{P}}_{\mathcal{D}^+}\left\{ \Psi \right\} (\frac{ \bm{\omega}_1}{2\pi})
&=&
 \frac{2 \pi C_o}{C_o^{1D}}{\mathcal{P}}_{\mathcal{D}}^{1D} \left\{
\Psi\right\}(\frac{\omega_1}{2\pi}) \\ \nonumber
&=&\lambda^{1D} \Psi(\frac{\omega_1}{2\pi}),
\end{eqnarray}
where we have used that $2 \pi C_o/C_{o}^{1D}=1$. Thus the solutions of (\ref{eigen1})
are given by
\begin{equation}
\Psi^{2D}_{n;l,m}(\bm{f})=\frac{1}{\sqrt{f}}\Psi^{(e)}_{n;\beta,\gamma}(f),\;f>0,
\end{equation}
where $\Psi^{(e)}_{n;\beta;\gamma}(\cdot)$ are the even Morse wavelets as defined
in one dimension \cite{Olhede2002}, $l=\beta-\frac{1}{2}$, $m=\gamma$, and
$n \in N$ enumerates the eigenvectors. 
The eigenvalues correspond to,
\begin{equation}
\lambda_{n,r}(C)=\lambda^{1D}_{n,r}
= \frac{\Gamma(r+n)}{\Gamma(n+1)\Gamma(r-1)}\int_0^{\frac{C-1}{C+1}}
x^n (1-x)^{r-2}\;dx,\end{equation}
and this yields the radial-spatial, radial-scale concentration of  $\mu_{\psi_{n}^{(e)}}
({\mathcal{D}(C)})=\lambda_{n,r}^2(C).$
It may seem surprising that, in two-dimensions, the same eigenvalues, and
thus concentration values,
are found as in the one-dimensional case, note, however, that this is only
derived for radial images, that are constrained
to the same behaviour in both spatial directions -- hence in essence we are
really only making a one dimensional compromise.
The {\em hypervolume} of ${\mathcal{D}}$ is directly related to $C,$ and
note that we may formulate the notion of bias of estimation of local properties
of the signal, or leakage,
in terms
of the eigenvalues as $1-\lambda_{n,r}^2(C)$ \cite{Olhede2002}.

\subsection{Non-Radial Localisation Operators}
The operator outlined above was constructed in the radial frequency domain
for explicitly radial images. Let us explicitly consider how a operator is constructed in two dimensions. Define
\begin{equation}
{\mathcal{D}}_2=\left\{\bm{\xi}:a^2+b^2+1 \le a 
C,\;0\le \theta \le 2\pi,\;\bm{b}\in{\mathbb{R}}^2
,\;a>0\right\}.
\end{equation}
Let us again consider a coherent state $\tilde{v}_{\bm{\xi}}^+(\bm{x})$ that
is the building
block of the projection operator, but let us now make this local to $\bm{\xi}.$
Define the Fourier transform of this coherent
state at $\bm{\xi}$ with $2\pi\bm{f}=\bm{\omega}$ as
\begin{eqnarray}
\tilde{V}_{\bm{\xi}}^+
(\bm{\omega})&=&
  a^{1/(2\gamma)+1/2}
 \omega^{-1/2+\gamma/2} V^{\beta,\gamma}(a^{1/\gamma}\omega)
e^{-j\bm{b} (\bm{r}_{-\theta}\bm{\omega}) \omega^{\gamma-1}}\nonumber \\
&=& \frac{2^{(r+1)/2}  \sqrt{\pi \gamma} }{\Gamma(r)} a^{(r+1)/2} \omega^{-1/2+\gamma/2}
\omega^{\beta} e^{-a \omega^{\gamma}-j\bm{b} (\bm{r}_{-\theta}\bm{\omega})
\omega^{\gamma-1}}.
\end{eqnarray}
The normalization of $\tilde{V}^+,$ for completeness, can be calculated
by:\\
$
\langle \tilde{v}_{\bm{\xi}}^+, \tilde{v}_{\bm{\xi}}^+\rangle=
 a^{1/\gamma+1} \frac{1}{(2\pi)^2}
\int_{0}^{\infty}\int_{0}^{2\pi} \omega^{-1+\gamma}V^{\beta,\gamma 2}(
a^{1/\gamma}\omega)
\omega d\omega d\phi  
= \frac{r}{2}.$
Therefore, the coherent states are of norm unity, if they are  multiplied
by $\sqrt{\frac{2}{r}}.$
We define the localisation operator for {\em any} function $g(\bm{x})\in
L^2({\mathbb{R}}^2)$ as
\begin{equation}
\tilde{{\mathcal{P}}}_{\mathcal{D}_2}\left\{g\right\}(\bm{x})
=\tilde{C}_o\int_{\mathcal{D}_2} \tilde{v}_{\bm{\xi}}^+(\bm{x})\langle 
\tilde{v}_{\bm{\xi}}^+, g\rangle\; \frac{da}{a^3}d^2\bm{b}d\theta.
\label{oppy2}
\end{equation}
Equation (\ref{oppy2}) gives an expression for the localisation of an arbitrary
function $g(\cdot)$ over region $\mathcal{D}_2.$ We can by calculating 
$\mu_{g}
({\mathcal{D}}_2)$ find the localisation of $g()\in L^2({\mathbb{R}}^2),$
to ${\mathcal{D}}_2.$
We shall take $\tilde{C}_o$ such that as $\left|{\mathcal{D}}_2\right|\longrightarrow
\infty$, 
$\tilde{{\mathcal{P}}}_{\mathcal{D}_2}\left\{g\right\}(x)\rightarrow g(x)$ for all
{\em radially symmetric} functions.
Note that if $g(\cdot)$ is a radially symmetric function, then
$
\langle 
\tilde{v}_{\bm{\xi}}^+, g\rangle
= \frac{1}{(2\pi)} 
\int_{0}^{\infty}  a^{1/(2\gamma)+1/2}
 \omega^{-1/2+\gamma/2} V^{\beta,\gamma}(a^{1/\gamma}
\omega) G(\frac{\omega}{2\pi}) J_0(\omega^{\gamma} b)\omega\;d\omega,
$
where $J_0(\cdot)$ is the zeroth Bessel function \cite{Abramowitz}.
Integration over $\phi$, for a radial $g(\cdot),$ thus removes the angular dependence on $\theta$ and $\phi_b,$
see \cite{McLachlan}.
In the frequency domain
\begin{eqnarray}
\nonumber
\tilde{{\mathcal{P}}}_{\mathcal{D}_2}\left\{G\right\}(\bm{\omega}_1)
&=& \tilde{C}_o \int_{\mathcal{D}} \tilde{V}_{\bm{\xi}}^+( \omega_1)\langle
\tilde{V}_{\bm{\xi}}^+, G\rangle\; \frac{da}{a^3}d^2\bm{b}d\theta\\
\nonumber
&=&2\pi \tilde{C}_o \int_{\omega_=0}^{\infty} \int \int_{a^2+b^2+1\le 2a C} a^{1/(2\gamma)+1/2}
\omega_1^{-1/2+\gamma/2} V^{\beta,\gamma}(a^{1/\gamma}\omega_1)  \omega_2^{1/2+\gamma/2}
a^{1/(2\gamma)+1/2}\\
&&\times V^{\beta,\gamma}(a^{1/\gamma}\omega_2)
 G(\frac{\omega_2}{2\pi}) J_0(b \omega_1^{\gamma}) J_0(b \omega_2^{\gamma}
 )\;\frac{da}{a^3} b db\;d\omega_2\\
&=&\tilde{C}_o \int_{\omega_2>0} \int \int_{a^2+b^2+1\le 2a C}
\kappa(\omega_1,\omega_2;a,b)\;\frac{da}{a^2}\;db\;d\omega_2.
\end{eqnarray}
Note that $ V^{\beta,\gamma} (\omega)\approx 0$ if $\omega<\omega_{\epsilon}$ and thus
we consider the term $J_0(\omega_l^{\gamma} b )$ for 
$\omega_l^{\gamma} b>>0,$ as the point $b=0$ has zero measure in the plane.
For $\left|z\right|>>0$ the asymptotic approximation to the zeroth Bessel function is  
\[J_0(z ) \approx \frac{\sqrt{2}}{\sqrt{\pi z}}\cos(z-\frac{\pi}{4}),
\]
{\em cf} \cite[p.~364,\;8.2.1]{Abramowitz}. We require $\left|\arg
z\right|<\pi$ however, this is not an issue as $b>0$ and $\omega_l>0.$
Thus for fixed non-zero $b$ for values such that the integrand is non-zero
\begin{eqnarray}
\nonumber
\kappa(\omega_1,\omega_2;a,b)&= & 2\pi a^{1/\gamma} 
 \omega_1^{-1/2+\gamma/2} V^{\beta,\gamma}(a^{1/\gamma}\omega_1)  \omega_2^{1/2+\gamma/2} V^{\beta,\gamma}(a^{1/\gamma}\omega_2) G(\frac{\omega_2}{2\pi})\\
\nonumber
&& J_0(b \omega_{1}^{\gamma}) J_0(b \omega_2^{\gamma} )b
\\
\nonumber
&\approx & 2\pi a^{1/\gamma}  \omega_1^{-1/2+\gamma/2}
V^{\beta,\gamma}(a^{1/\gamma}\omega_1)  \omega_2^{1/2+\gamma/2}
V^{\beta,\gamma}(a^{1/\gamma}(\omega_2) G(\frac{\omega_2}{2\pi}) \\
&& b \sqrt{\frac{2}{\pi \omega_2^{\gamma} b}}\cos(\omega_2^{\gamma} b-\frac{\pi}{4})
\sqrt{\frac{2}{\pi \omega_1^{\gamma} b}}\cos(\omega_1^{\gamma} b-\frac{\pi}{4})
\nonumber
\\
\nonumber
&=&2 a^{1/\gamma} 
V^{\beta,\gamma}(a^{1/\gamma}\omega_1) 
\sqrt{\omega_2/\omega_1}
V^{\beta,\gamma}(a^{1/\gamma}\omega_2) G(\frac{\omega_2}{2\pi})\cos((\omega_1^{\gamma}- \omega_2^{\gamma})b)
\\
&&+
2 a^{1/\gamma} 
V^{\beta,\gamma}(a^{1/\gamma}\omega_1) 
\sqrt{\omega_2/\omega_1}
V^{\beta,\gamma}(a^{1/\gamma}(\omega_2)) G(\frac{\omega_2}{2\pi}) \cos((\omega_1^{\gamma}+ \omega_2^{\gamma})b-\frac{\pi}{2}).
\label{approx}
\end{eqnarray}
The integration is over $\omega_2>0$ and also $\omega_1>0$, hence the first term of (\ref{approx}) dominates over the second term. The integrand of $\kappa(\omega_1,\omega_2;a,b)$ can be replaced by
$\kappa_2(\omega_1,\omega_2;a,b)$ defined in (\ref{kappa2}) of the previous section.
Thus we approximate the operator acting on $G(\cdot)$ in the frequency domain
via
\begin{eqnarray}
\tilde{{\mathcal{P}}}_{\mathcal{D}_2}\left\{G\right\}(\bm{\omega}_1)&\approx&
\tilde{C}_o \int_{\omega_2=0}^{\infty} \int \int_{a^2+b^2+1\le 2a C}\kappa_2(\omega_1,\omega_2;a,b)
\;\frac{da}{a^2}\;db\;d\omega_2\equiv
{\mathcal{P}}_{\mathcal{D}}\left\{G\right\}(\bm{\omega}_1),
\end{eqnarray}
and $\tilde{C}_o$ can be found from $C_o.$
Defining $\tilde{{\mathcal{P}}}_{\mathcal{D}_2}$ allows for the consideration
of the localisation of an arbitrary $L^2({\mathbb{R}}^2)$ function. The radial
eigenfunctions of $\tilde{{\mathcal{P}}}_{\mathcal{D}_2}$ are approximately
those of ${\mathcal{P}}_{\mathcal{D}},$ where the derivation of the approximation
shows the reasoning behind the definition of ${\mathcal{P}}_{\mathcal{D}}.$
Finally $\tilde{{\mathcal{P}}}_{\mathcal{D}_2}$ can be generalised to an
arbitrary localisation by removing the constraint of radial symmetry in ${\mathcal{D}_2}$
and on $g(\cdot).$
%
 
\subsection{Isotropic Wavelet Definition}
The even multiple Morse wavelets are defined in one dimension in the Fourier domain
for fixed $n=0,1,2\dots,\;\beta\ge 1,\;\gamma> (\beta-1)/2$ and denoted
by $\Psi^{(e)}_{n;\beta,\gamma}\left(f\right)$
in terms of  $A_{n;\beta,\gamma}=\sqrt{\pi \gamma 2^r\Gamma(n+1)/\Gamma(n+r)}$
\cite{Olhede2002}.
We define the isotropic two dimensional wavelets as the eigenfunctions of
${\mathcal{P}}_{\mathcal{D}},$ for
fixed $n,l,m$ in terms of $f=\|\bm{f}\|$ via
\begin{equation}
\Psi^{(e)}_{n;l,m}\left(\bm{f}\right)
= \frac{A_{n;l,m}}{\sqrt{\pi}} 
(2\pi f)^l e^{-(2 \pi f)^m}L_{n}^{c^{\prime}_{l,m}}\left(2 (2\pi f)^{m}\right),
\label{ftpsi}
\end{equation}
with $c^{\prime}_{l,m}=(2l+2)/m-1,$  to attain the correct normalisation
over $L^2({\mathbb{R}}^2)$ and where $L_{n}^{c}\left(\cdot\right)$ are generalized
Laguerre polynomials.
The spatial domain wavelets with $x=\|\bm{x}\|$ are found via the inverse Fourier transform for radial images, 
$\psi^{(e)}_n\left(\bm{x}\right)=2\pi \int_{0}^{\infty} \Psi^{(e)} \left(f\right)
J_0\left(2\pi f x\right) f df.$
We plot the wavelets for $l=8$ and $m=3$ with $n=0,1$ in the spatial domain,
see Figure \ref{fig:partition}. Their radially symmetric oscillatory structure is very
clear. They are optimally concentrated in a radial structure centred at
the origin.
We plot the modulus of the same function in the spatial frequency
domain in the same plot. These functions are band-pass filters that
are non-zero for a range of frequencies centred at the same distance from
the origin in the frequency domain. The trough in the $n=1$ follows as the
first two wavelets are orthogonal, and we see that the sum of the moduli
will be large in the same ring-shape structure.
To further characterise the Morse wavelets define the radial frequency that maximises the isotropic
Morse wavelets as
\begin{equation}
\label{fmaxxy}
f_{\mathrm{max}}^{(n)}=\arg_{f>0}\max \left|\Psi_{n}^{(e)}
\right|^2,\;n=0,\dots,N-1.
\end{equation}
The magnitude square of the Fourier transform of analysis wavelet 
$\psi_{\bm{\xi},n}^{(e)}(\cdot)$
will have a maximum at frequency $f_{\mathrm{max}}^{(n)}/a$ and is unaffected by both $\bm{b}$
and the rotation. 

\section{Monogenic Images}
In one dimension the analytic signal is used to unambiguously define the phase and 
amplitude of a given real signal and using an analytic analysis wavelet will allow for the definition
of a local magnitude and phase at each time and scale  point -- a necessity
for the analysis of multi-component images. 
The analytic signal is constructed
in one dimension by removing all negative frequencies in the signal, and then inverting
the Fourier transform -- any real signal $u(t)$ is complemented by its
Hilbert transform $v(t)={\mathcal{H}}\left\{u\right\}(t),$ and the analytic
signal corresponds to $u^+(t)=u(t)+iv(t).$ If an oscillation is persistent
at a particular range of times,
then it will be heavily weighted in the Fourier domain, and the analytic
signal will
approximately take the form of a complex exponential. 
\\
The correct extension of the analytic signal to two dimensions has been the
subject of much debate -- of particular note are perhaps the single orthant
image of \cite{Hahn1992}, the hypercomplex signal of \cite{vonBulow} and the monogenic signal of \cite{FelsbergSommer}. 
Following Felsberg and Sommer \cite{FelsbergSommer} we define the Riesz transform of 
an image $u(\bm{x})$ as 
\begin{equation}
\label{rieszy}
{\mathcal{R}}u(\bm{x})=\bm{i}{\mathcal{R}}_1u(\bm{x})+\bm{j}{\mathcal{R}}_2u(\bm{x})=
\bm{i} v^{(1)}(\bm{x})+\bm{j} v^{(2)}(\bm{x}).
\end{equation}
The Fourier transforms of these two objects are
\[{\mathcal{F}}\left\{{\mathcal{R}}_1u(\bm{x})\right\}=
-\bm{j}\cos\left(\phi\right)U(\bm{f}),\;{\mathcal{F}}\left\{{\mathcal{R}}_2u(\bm{x})\right\}=
-\bm{j}\sin\left(\phi\right)U(\bm{f}).\]
Define the monogenic image \cite{FelsbergSommer}
as $u^{(+)}(\bm{x})= u(\bm{x})+{\mathcal{R}}u(\bm{x}).$
This is a quaternion valued object, and relations between the components
of the quaternion are interpretable in
terms of orientation and phase, for oscillatory images, as will be demonstrated.

Consider an oscillation in two dimensions, corresponding to repeating maxima
spaced $1/f_0$ apart in 
orientation $\bm{n}=\left[\cos(\eta)\;\sin(\eta)\right]^T.$
This corresponds to the simplest oscillatory image
\begin{equation}
\label{ossy1}
g_1(\bm{x})=a_1\cos\left(2\pi f_0 \bm{x} \cdot \bm{n}+\theta_s\right),
\end{equation}
where $a_1,$ $f_0$ and $\theta_s$ are constant scalars, whilst $\bm{n}$ is
a constant unit length
vector. Note that in variable $y=\bm{x} \cdot \bm{n},$ $g_1(\bm{x})$
solves the one dimensional harmonic oscillator equation of
$
\left[\frac{\partial^2}{\partial y^2}+
\left(2\pi f_0\right)^2\right]g_1(y)=0.
$
To characterise a given image, we wish to determine $a_1,$ $f_0,$ $\theta_s$
and $\bm{n},$ from the image.\\
We can calculate the monogenic extension of (\ref{ossy1}), and
obtain, with $\bm{e}_{\eta}=\bm{i}\cos(\eta)+\bm{j}\sin(\eta)$, \cite{Felsbergtech}
\begin{eqnarray}
\label{monny}
g_1^{(+)}(\bm{x})&=& g(\bm{x})+\bm{i}g^{(1)}(\bm{x})+\bm{j}g^{(2)}(\bm{x}) \nonumber \\ 
               &=& a_1\left[\cos\left(2\pi f_0 \bm{x} \cdot \bm{n}+\theta_s\right)
+\bm{e}_{\eta}\sin\left(2\pi f_0 \bm{x} \cdot \bm{n}+\theta_s\right)\right].
\end{eqnarray}
In quaternionic polar coordinates, 
\begin{equation}
g^{(+)}_1({\bf x}) = a_1 e^{2\pi \bm{e}_{\eta} (f_0 y+\theta_s)},
\end{equation}
and we may determine
$
a_1= \sqrt{ g^2(\bm{x})+ g^{(1)2}(\bm{x}) + g^{(2)2}(\bm{x})   },$ 
$\eta = \tan^{-1}\left(\frac{g^{(2)}(\bm{x})}{g^{(1)}(\bm{x})}
\right),$\\
$ f_0 y+\theta_s = \frac{1}{2\pi}
\tan^{-1}\left(\frac{{\mathrm{sgn}}\left(g^{(1)}(\bm{x}) \right)
\sqrt{g^{(1)2}(\bm{x})+g^{(2)2}(\bm{x})}}{g(\bm{x})}\right).
$
We restrict $-\frac{\pi}{2}<\eta<\frac{\pi}{2},$ and
$-\frac{1}{2} < f_0 y+\theta_s \le \frac{1}{2}.$ 
In terms of $y$ we hence have a frequency domain description that perfectly
mirrors the one dimensional theory, and retrieve the properties
of the image from its monogenic extension.
Naturally, in most real applications 
 perfectly oscillatory images are not encountered, and more general
models must be considered.

\section{AM/FM/OM Images}
As in  \cite{Gonnet} and \cite{Havlicek}, we shall consider images that locally
may be approximated as a sum of sinusoidal components. This model is of some
importance in machine vision, and can be applied to granular flow and general
oriented patterns \cite{Bovik}.
Assume that
\begin{equation}
c(\bm{x})=\sum_{l=1}^L c_l(\bm{x}),\;c_{l}(\bm{x})=a_l(\bm{x})\cos(2\pi
\phi_l(\bm{x})),
\label{amfmom}
\end{equation}
where $\phi_l(\bm{x})= \varphi_l(\bm{n}_l(\bm{x})\cdot\bm{x}),$ 
with the added constraint that the unit vector $\bm{n}_l(\bm{x})$ is varying
slowly, in comparison
with $\bm{x}$, across the spatial period. We introduce this extra notation
so that we may characterise images that can be considered as approximately
sinusoidal in variable $y_l(\bm{x})=\bm{n}_l(\bm{x})\cdot\bm{x}.$ 
We refer to $\bm{n}_l(\bm{x})$
as the orientation modulation (OM), whilst $\varphi_l(\cdot)$ is the phase
modulation of component $l,$ and plays the same role as the phase/frequency modulation (FM)
of a one dimensional signal.
Then we find that for $\bm{x}=\bm{x}_o+\bm{\delta} \bm{x},$ and 
$\varphi_l^{\prime}(y_l)=
d/dy_{l} \left[ \varphi_{l}(y_l) \right],$ that
\begin{equation}
\phi_l(\bm{x})=\phi_l(\bm{x}_o)+
\varphi_l^{\prime}(\bm{n}_l (\bm{x}_o)\cdot\bm{x}_o)
\bm{n}_l (\bm{x}_o)\cdot\left(\bm{x}-\bm{x}_o\right)+
O\left(\left| \bm{\delta} \bm{x} \right|^{2} \right).
\end{equation}
In the following we will use the shorthand $\phi_{l}^{\prime}(\bm{x} )
= \varphi_l^{\prime}( \bm{n}_l(\bm{x})\cdot  \bm{x})$. 
We additionally assume $a_{l}(\bm{x})$
varies slowly in
comparison to the
cosine term, and this corresponds to the amplitude modulation (AM). The monogenic
version of the $l$th component is
\begin{eqnarray}
\label{monogenicloc}
c^+_l(\bm{x})&=&a_l(\bm{x}_o) e^{2\pi \bm{e}_{l}(\bm{x}_o)
\left(\phi_l(\bm{x}_o)+
\phi_l^{\prime}(\bm{x}_o)
\bm{n}_l (\bm{x}_o)\cdot \left(\bm{x}-\bm{x}_o\right)+O\left(\left|\bm{x}-\bm{x}_o\right|^2 \right) \right)}\\
\nonumber
&=&a_l(\bm{x}_o)
\left[
\cos\left(2\pi \phi_l( \bm{x}_o)+2\pi
\phi_l^{\prime}(\bm{x}_o) 
 \bm{n}_l(\bm{x}_o) \cdot\left(\bm{x}-\bm{x}_o\right) \right)
\right.\\
&&
\left.
+{\bm{e}}_{l}({\bm{x}}_o)
\sin\left(2\pi\phi_l(\bm{x}_o)+2\pi
\phi_l^{\prime}(\bm{x}_o)
\bm{n}_l (\bm{x}_o)\cdot\left(\bm{x}-\bm{x}_o\right) \right)
\right]+o(\left|\bm{x}-\bm{x}_o\right|).
\nonumber
\end{eqnarray}
{\em cf} equation (\ref{monny}).
\begin{figure*}[t]
\centerline{
\includegraphics[height=1.75in,width=1.75in]{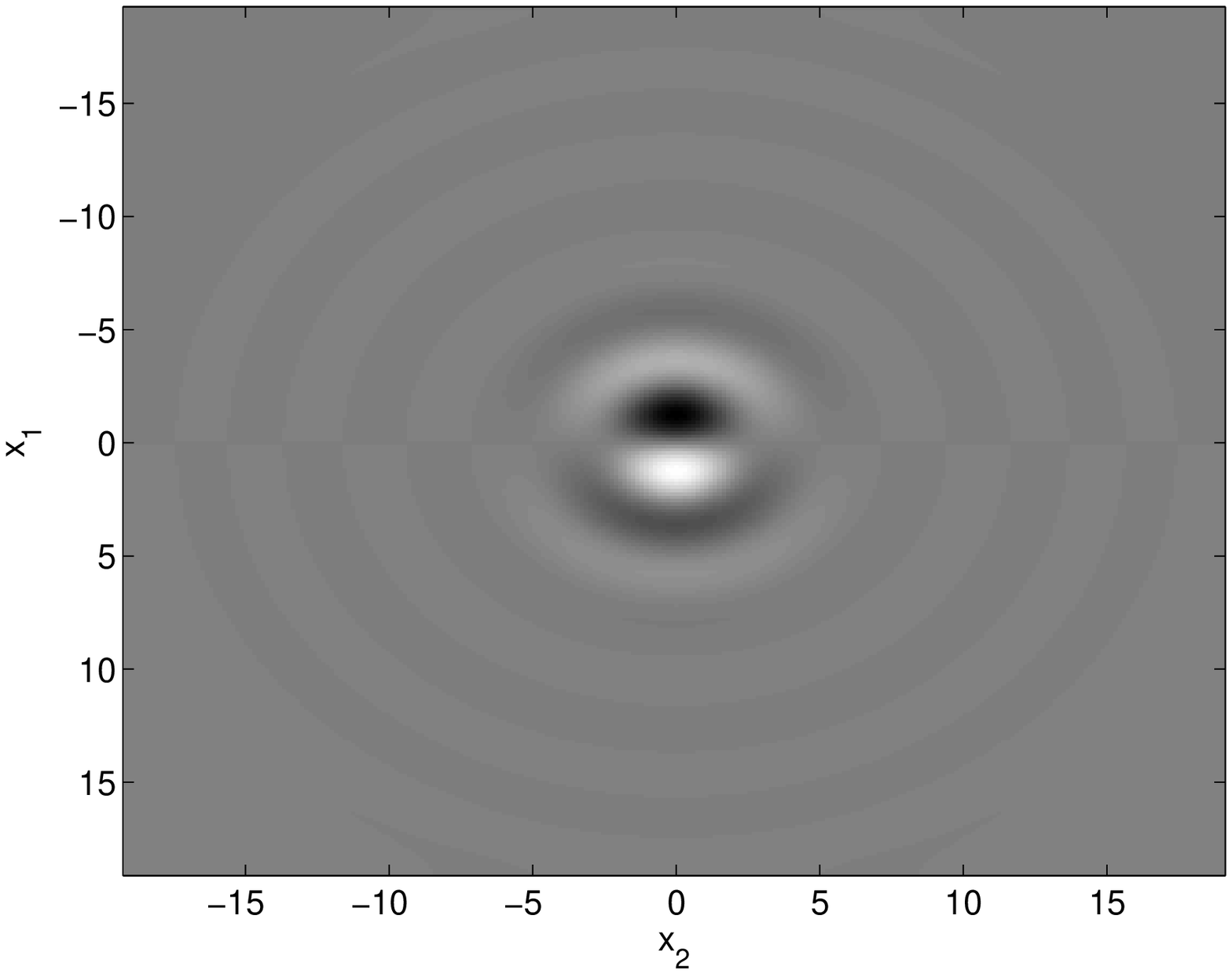}
\includegraphics[height=1.75in,width=1.75in]{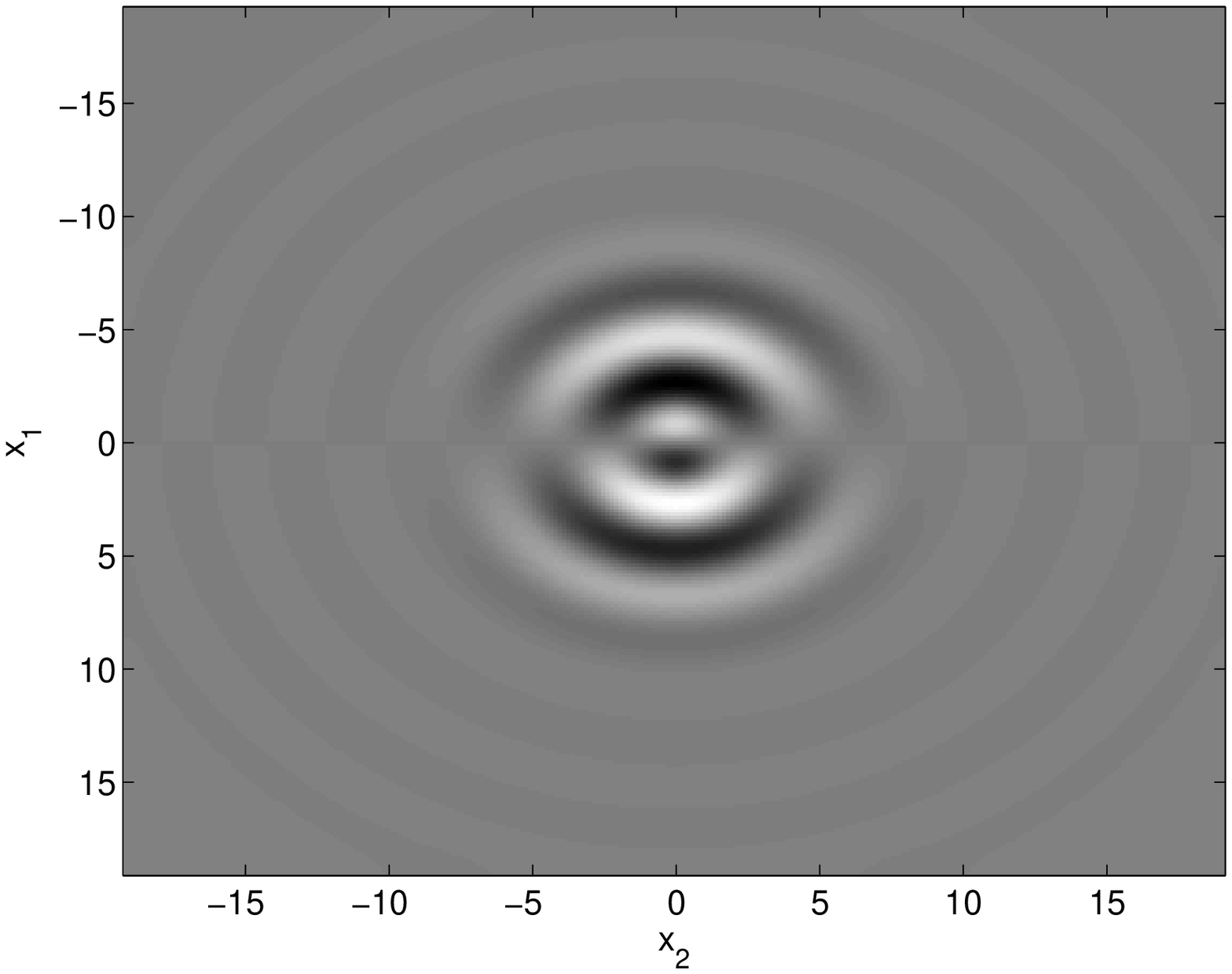}
\includegraphics[height=1.75in,width=1.75in]{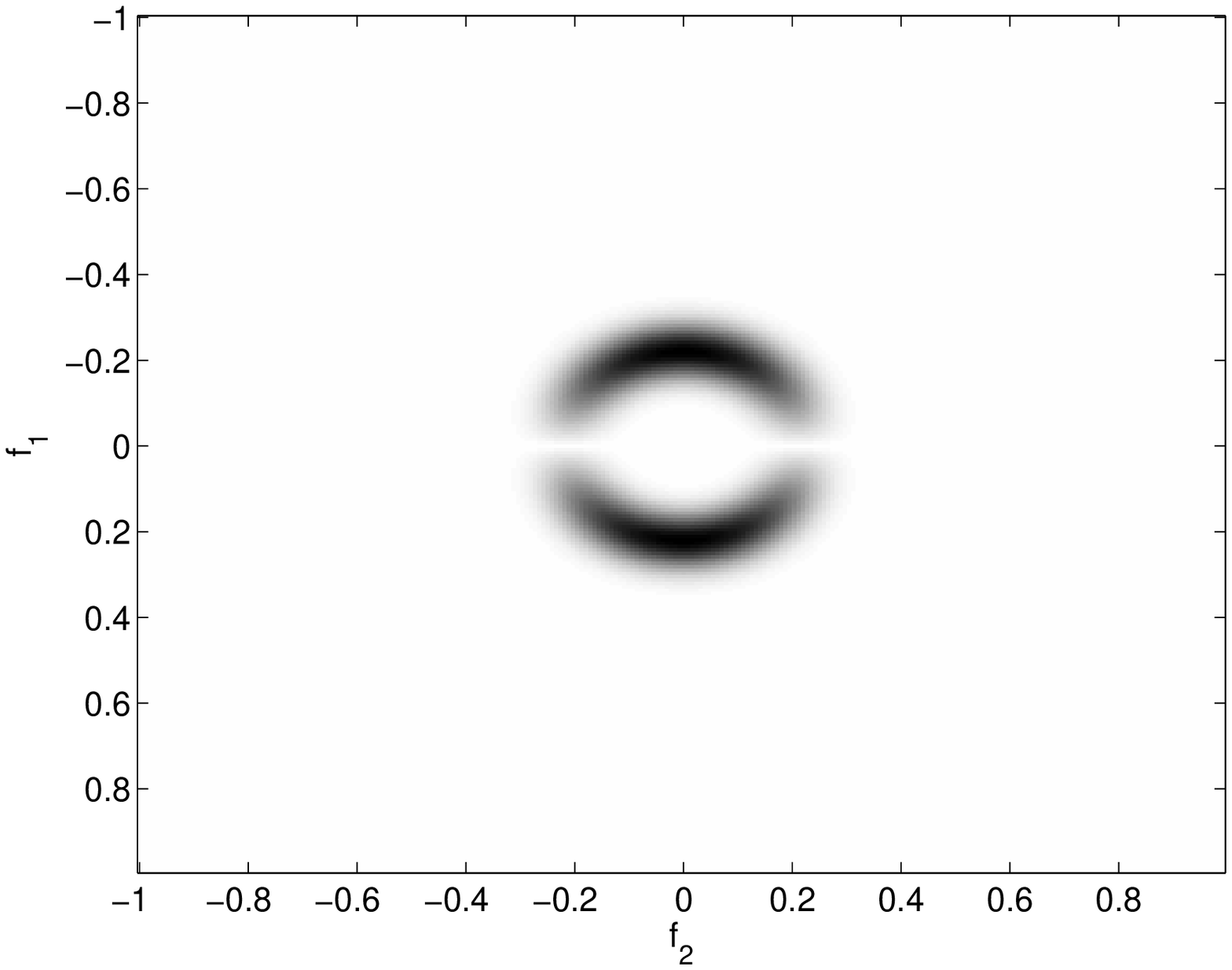}
\includegraphics[height=1.75in,width=1.75in]{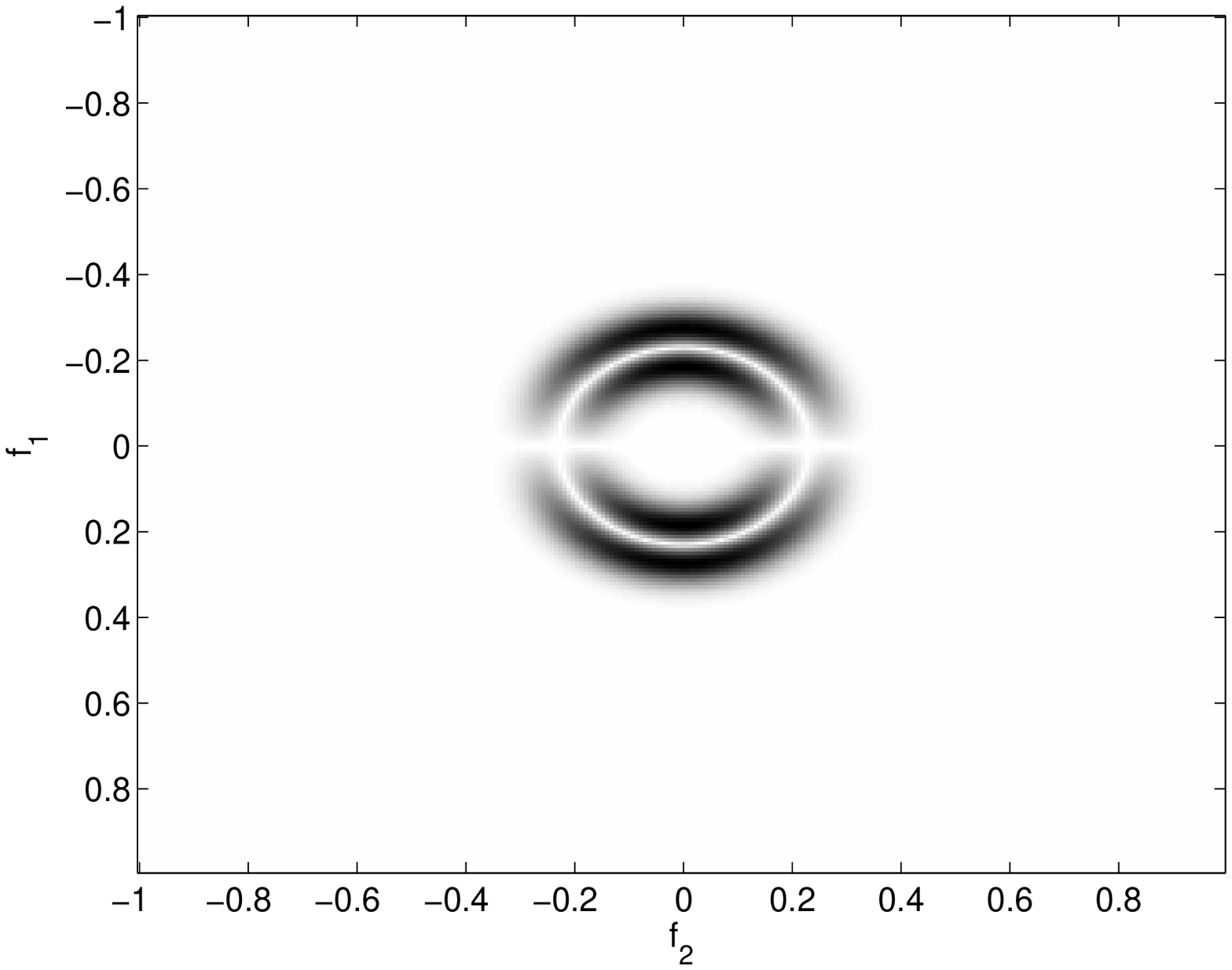}
}
\caption{\label{fig:partition5} 
The $x_1$ Riesz transform Morse wavelets in the spatial
domain for $(l,m)=(8,3)$ and $n=0$ (far left), $n=1$ (second left).
The modulus of the $x_1$ Riesz transform
Morse wavelets in the spatial frequency domain for $(l,m)=(8,3)$ and $n=0,$
(second from right)
$n=1$ (far right).
}
\end{figure*}
%
%
Assuming the orientation $\bm{e}_l(\bm{x}_o)$ to be stable across values
of $\bm{x}$ for which $a_l(\bm{x})$ is non-zero,
we shall perform the Fourier transform of $c_{l}^{(+)}(\bm{x})$ in terms of
the unit quaternion $\bm{e}_l(\bm{x}_o)$
instead of $\bm{j}$.
 Note that de Moivre's theorem is still valid for any
unit quaternion, and so this directional Fourier transform can be interpreted
just like
the regular Fourier transform in terms of oscillatory components. 
The directional Fourier transform is then given by
\begin{equation}
\label{quatfft}
C^+_{\bm{e}_l,l}(\bm{f}) =\int \int a_l( \bm{x}_o )
 \; e^{2\pi \bm{e}_l (\bm{x}_o) \left( \phi_l(\bm{x}_o)+
 \phi_l^{\prime}(\bm{x}_o) \bm{n}_l(\bm{x}_o) (\bm{x}-\bm{x}_o) \right)}
 e^{-2 \pi \bm{e}_l(\bm{x}_o) \bm{f} \bm{x}}\;d^2\bm{x}.
\end{equation}
We apply the stationary phase approximation to this integral \cite{Gonnet},
as the unit quaternion will be acting like $\bm{j}$ itself. Thus, under
the assumption that $\Omega(\bm{x},\bm{f})=\phi_l(\bm{x}_o)+
\phi_l^{\prime}(\bm{x}_o)\bm{n}_l (\bm{x}_o) \cdot (\bm{x}-\bm{x}_o)
- \bm{f} \bm{x}$ has the unique stationary
point $\bm{x}_o$, and there is quadratic behaviour around this point, i.e.
$\Omega(\bm{x},\bm{f})=\frac{1}{2}\left[(\bm{x}-\bm{x}_0)\bm{H}(\bm{x}-\bm{x}_0)
\right],$
where $\bm{H}$ is the Hessian matrix of $\Omega(\bm{x},\bm{f}),$ 
we find at 
$\bm{f}(\bm{x}_o) = \phi_l^{\prime}(\bm{x}_o) \bm{n}_l (\bm{x}_o),$
the integral provides the only non-null contribution of 
$C^+(\bm{f})\approx \frac{2\pi}{\sqrt{\left|\bm{H}\right|}}
a_l(\bm{x}_o) e^{2\pi \bm{e}_l( \bm{x}_o) \Omega(\bm{x}_o, \bm{f})}$.
This provides a local frequency description of component $l$ at $\bm{x}_o.$
The instantaneous frequency should be interpreted in terms of
the local frequency of the oscillations, $\left|
\phi_l^{\prime}(\bm{x}_o)\right|,$ and the orientation of these oscillations,
 $\bm{n}_l(\bm{x}_o).$ The sign of $\phi_l^{\prime}(\bm{x}_o)$ is taken
so that the orientation angle is restricted from $-\frac{\pi}{2}$ to
$\frac{\pi}{2}.$  
Finally the local magnitude $
a_l(\bm{x}_o)$ has the interpretation of local energy presence.
If the image is actually a sum of several AM/FM/OM terms, {\em i.e.} the
signal is multi-component, as is most often
the case of many observed images, we will not
be able to use this description directly, as we cannot separate out the $L$
components. The problem with multi-component signals
in one dimension is much documented \cite{Boashash1992}, if not fully resolved,
but generally calls for localised methods. We shall
thus construct local monogenic descriptions of images. These descriptions extract well-behaved orientation, phase and
amplitude functions locally, and also have excellent statistical properties.

%

\section{\label{monor} Orientation \& Monogenic Wavelets}
\subsection{Definition}
In one dimension the analytic Morse wavelets can be constructed from the
even
Morse wavelets, by adding $\bm{i}$ times the Hilbert transform of the original
even function, to the even wavelet. The even Morse wavelets
are invariant to sign changes, or direction, whilst
the odd Morse wavelets, are naturally odd functions. 
We construct
a monogenic version of isotropic wavelets in two dimensions, based on the
real isotropic Morse wavelets, using the Riesz transform.
The monogenic wavelet transform 
locally defines a phase and an orientation for a real image at each spatial
and spatial scale point, similarly to the monogenic signal.
The monogenic wavelets
are most easily represented as quaternion-valued functions defined for each
$n$ via
\begin{equation}
\psi^{(+)}_n (\bm{x})  = \psi^{(e)}_n\left(\bm{x}\right)+\bm{i}
\psi^{(1)}_n\left(\bm{x}\right)+\bm{j}
\psi^{(2)}_n\left(\bm{x}\right)=\psi^{(e)}_n\left(\bm{x}\right)+\psi^{(q)}_n\left(\bm{x}\right),
\end{equation}
where the Fourier transform of the real part of the monogenic wavelet is given by (\ref{ftpsi}) 
and the additional two real functions are defined in the Fourier domain
as Riesz transforms of the isotropic wavelet function
via
\begin{equation}
\Psi^{(s)}_{n,l,m} \left(\bm{f}\right)
= -\bm{j}\frac{A_{n;l,m}}{\sqrt{\pi}} 
\frac{f_s}{f} (2\pi f)^l e^{-(2 \pi f)^m}L_{n}^{c^{\prime}_{l,m}}\left(2
(2\pi f)^{m}\right),\;s=1,2.
\end{equation}
Note that we consider 
$
\psi^{(q)}_n\left(\bm{x}\right)=\bm{i}
\psi^{(1)}_n\left(\bm{x}\right)+\bm{j}
\psi^{(2)}_n\left(\bm{x}\right),
$
as a single object, and this has the same norm as $\psi^{(e)}_n\left(\bm{x}\right)$.
In subsequent analysis we fix $(l,m),$ and henceforth suppress their value,
for notational convenience. 
For the case $(l,m)=(8,3)$ plots of the Riesz transforms in the $x_1$ direction
of the isotropic wavelets
$n=0,1$ are given in the spatial domain (see Figure \ref{fig:partition5})
as well as the spatial frequency domain, where
their modulus is plotted.
Note that the real component is an isotropic wavelet (like the even wavelet
in one dimension) and the two components are odd in the $x_1$ and $x_2$ direction
respectively.
We define the translated, rotated and dilated wavelet 
$\psi^{(+)}_{\bm{\xi},n} (\bm{x}),$ as the appropriate sum of translating,
rotating and dilating its real valued components.
The 
continuous wavelet transform of an image $g(\cdot)$ with respect to either
the components of, or with respect to the full quaternionic wavelet, is defined
as
$w_n^{(\cdot)}(\bm{\xi};g)= \int d^{2} \bm{x} g(\bm{x}) 
\psi^{(\cdot)*}_{\bm{\xi},n}\left(\bm{x}\right)$.
The associated scalogram is $S_n^{(\cdot)}(\bm{\xi};g)= |w_n^{(\cdot)}(\bm{\xi};g)|^2.$
The wavelet transform of image $g(\cdot)$
is then given via either the spatial domain, or spatial frequency domain,
in terms of 
$\bm{\zeta}=\left[a,\theta,\bm{f}_{\bm{b}}\right]^T$, where $\bm{f}_{\bm{b}}$
is the Fourier variable for $\bm{b}$ via
\begin{eqnarray}
w^{(+)}_n\left(\bm{\xi};g\right)&=&  w^{(e)}_n\left(\bm{\xi};g\right)-\bm{i}
w^{(1)}_n\left(\bm{\xi};g\right)-\bm{j} 
 w^{(2)}_n\left(\bm{\xi};g\right), \\
W^{(+)}_n\left(\bm{\zeta};g\right)&=&  
\left[1-\bm{k}\cos(\phi-\theta)+\sin(\phi-\theta) \right]
W^{(e)}_n\left(\bm{\zeta};g\right).
\end{eqnarray}
Thus the filtering carried out in the Fourier domain, can be understood,
by looking at the wavelet function in the Fourier domain:
\begin{equation}
\Psi^{(+)}_{\bm{\xi},n}(\bm{f})=\left[
1-\bm{k}\cos(\phi-\theta)+\sin(\phi-\theta)
\right] \Psi^{(e)}_{\bm{\xi},n}\left(\bm{f}\right), 
\end{equation}
with modulus 
\begin{equation}
\left|\Psi^{(+)}_{\bm{\xi},n}(\bm{f})\right|^2=2\left(
1+\sin(\phi-\theta)
\right)\left|\Psi^{(e)}_{\bm{\xi},n}\left(\bm{f}\right)\right|^2.
\end{equation}
Hence in terms of $\bm{b}$ and $a$ the monogenic wavelet is filtering the
image identically to the isotropic real wavelet, whilst the term $\sin(\phi-\theta)$
is positioning the image, in orientation, in relation to the wavelet. This
representation clarifies
that the rotation by angle $\theta$ of the wavelet function, as repositioning
the axis of analysis, by a rotation of $\theta.$ Note that the modulus of the real isotropic wavelet is invariant with respect to $\theta.$
\\
Finally consider the joint structure of the $N$ wavelets. 
We obtain that (see Appendix A)
\begin{equation}
\begin{array}{llllll}
\langle \psi_{n_1 }^{(e)}, \psi_{n_2}^{(e)}\rangle &=& \delta_{n_1,n_2}, & \langle \psi_{n_1 }^{(l)}, \psi_{n_2}^{(l)} \rangle &=& \frac{1}{2}\delta_{n_1,n_2}, \\
\langle \psi_{n_1 }^{(1)}, \psi_{n_2}^{(2)}\rangle &=& 0, & \langle \psi_{n_1 }^{(e)}, \psi_{n_2}^{(l)}\rangle&=&0,  \ \ l=1,2.
\end{array}
\end{equation}
 The multiple Morse wavelets thus
form an orthogonal system, and this will have
implications for their usage when performing estimation of local characteristics
of real images. 
The total energy of the image using the $n$th wavelet is given by
\begin{eqnarray}
\nonumber
S^{(+)}_n\left(\bm{\xi};g\right)&=&S^{(e)}_n\left(\bm{\xi};g\right)+
S^{(1)}_n\left(\bm{\xi};g\right)+
S^{(2)}_n\left(\bm{\xi};g\right)\\
&=&S^{(e)}_n\left(\bm{\xi};g\right)+
S^{(q)}_n\left(\bm{\xi};g\right).
\end{eqnarray}
Also note that,
as $\psi^{(e)}(\cdot)$ is radially symmetric,
\begin{eqnarray}
\nonumber
w_n^{(1)}(\bm{\xi};g)&=&\cos(\theta)w_n^{(1)}(\bm{\xi}_0;g)+
\sin(\theta)w_n^{(2)}(\bm{\xi}_0;g),\\
w_n^{(2)}(\bm{\xi};g)&=&-\sin(\theta)w_n^{(1)}(\bm{\xi}_0;g)+
\cos(\theta)w_n^{(2)}(\bm{\xi}_0;g),
\label{rotenergy}
\end{eqnarray}
where $\bm{\xi}_0=(a,0,\bm{b})$. Thus the wavelet transform needs only be calculated for one orientation, and can then be formed
for any orientation by judicious recombination. We consider now the analysis
of typical image features with the monogenic Morse wavelets.

\section{The Monogenic Wavelet Transform of Discontinuities}
We consider both point and line discontinuities.
The idealised version of a point discontinuity
at $\bm{x}_{0}$ corresponds to:
$g_{s,1}(\bm{x})=a_1(\bm{x})\delta(\bm{x}-\bm{x}_0),$
where $a_1(\bm{x})$ is assumed to be a well-behaved function at point $\bm{x}=\bm{x}_0.$
This singularity is characterised by $\bm{x}_0,$ its location, and 
$a_1(\bm{x}_0),$ the amplitude of the location. 
The wavelet transform of this object is
$
w^{(+)}(\bm{\xi};g)=a_1(\bm{x}_0)\psi^{(+)*}_{\bm{\xi}}\left(\bm{r}_{-\theta}
\left(\bm{x}_0-\bm{b})/a
\right)\right)/a.
$
This is clearly maximum near $\bm{b}=\bm{x}_0$ where it has a modulus 
square of 
$a_1^2(\bm{x}_0)\left|\psi^{(+)}_{\bm{\xi}}\left(\bm{0}\right) \right|^2/a^2,
$
and hence point singularities can be located by finding maxima in $\bm{b}.$
Furthermore $\left|w^{(+)}(\bm{\xi};g)\right|^2$ has no dependence on $\theta,$
and as the magnitude of the wavelet at the origin
is known, $a_1(\cdot)$ can be determined. A one-dimensional singularity is
modelled as
\begin{equation}
g_{s,2}(\bm{x})=a_2(x_1)\delta\left(\cos(\theta_2)x_1+\sin(\theta_2)x_2-c\right).
\label{dis}
\end{equation}
The line $x_2=c$ ($\theta_2 = \pi/2$) modulated by the value of $a_2(x_1)$
is permitted, however we
do not permit the line $x_1=c$ ($\theta_2=0$), as this would lead to an image
of infinite energy. We will
anticipate a further constraint on $\theta_2$, that is $-\frac{\pi}{4}
<\theta_2 \le \frac{3\pi}{4}$, which we will explain at the end of this section.
Assume that $a_2(x_1)$ is a symmetric function around a maximum at $x_1=x_{1,\max}.$
We characterise the structure of $g_{s,2}(\cdot),$ using the wavelet
transform. The wavelet transform using the isotropic wavelet only, noting that $\psi^{(e)}_n(\bm{x})=
\psi^{(e)}_{n,r}(x),$ (where $\psi^{(e)}_{n,r}(x)\approx 0\;\forall x>x_r$)
is \\
$
w_n^{(e)}(\bm{\xi},g_{s,2})=\int_{-\infty}^{\infty} a_2(x_1)/a 
\psi^{(e)}_{n,r}\left(a\sqrt{(x_1-b_1)^2+(c\csc(\theta_2)-\cot(\theta_2)x_1-b_2)^2
}\right)\;dx_1.$ This will be large for values of $\bm{b}$ such that 
$\sqrt{(x_{1,\max}-b_1)^2+(c\csc(\theta_2)-\cot(\theta_2)x_{1,\max}-b_2)^2}<x_r/a.$
Similar results hold for $w_n^{(e)}(\bm{\xi},g_{s,2}),$ as the Riesz transforms
roughly have the same spatial support as the original wavelet. Hence we identify
the location of the singularity for any fixed value $a$ as $\bm{b}=(x_{1,\max}, c\csc(\theta_2)-\cot(\theta_2)x_{1,\max} ),$ as maxima in
the modulus of the wavelet transform using the monogenic wavelet. The orientation $\theta_2$ will thus
visually be apparent from
the $\bm{b}$ plane but can also be characterised at a {\em fixed} point $\bm{\xi}.$
The Fourier transform of (\ref{dis}) is 
\begin{eqnarray}
\nonumber
G_{s,2}(\bm{f})
&=&\frac{A_2(f_1-\cot(\theta_2)f_2)}{\left|\sin(\theta_2)\right|} 
e^{-\bm{j}2\pi (f_1\cos(\theta_2)+f_2\sin(\theta_2))c}
e^{-\bm{j}2\pi (-\sin(\theta_2)f_1+\cos(\theta_2)f_2)c\cot(\theta_2)}
\end{eqnarray}
where $A_2(\cdot)$ is the Fourier Transform of $a_2(\cdot)$. The wavelet transform of this image is
$W_{n}^{(l)}(\bm{\zeta}_0;g_{s,2})= \bm{j} \frac{ f_{b,l} }{f_b} a\Psi^{(e)}_n(a
f_b)
G_{s,2}(\bm{f}),$
where $\bm{\zeta}_{0}=\left[ a, 0, \bm{f}_{\bm{b}} \right]$. 
Thus we find that the Fourier transform of the rotated wavelet is given by
$
W_{n}^{(1)}(\bm{\zeta};g_{s,2})= \bm{j} 
\frac{\cos(\theta)f_{b,1}+\sin(\theta)f_{b,2}}{f_b}a\Psi^{(e)}_n(a f_b)
G_{s,2}(\bm{f}),$
with the obvious extension for $W_{n}^{(2)}(\bm{\zeta};g_{s,2}).$
Now we choose to evaluate the wavelet transform with $\theta_2=\theta,$
and define $\bm{\xi}_2=\left[ a, \theta_2, \bm{b}\right].$
\begin{eqnarray*}
\nonumber
w_{n}^{(1)}(\bm{\xi}_2;g_{s,2})
&=&\frac{1}{\left|\sin(\theta)\right|} \int_{-\infty}^{\infty}
\bm{j} \frac{f_{\theta,b,1}}{f_{\theta,b}}
e^{\bm{j}2\pi f_{\theta,b,1}\left[(\bm{r}_{-\theta} \bm{b})_1-c\right]} \\
&&
 \int_{-\infty}^{\infty}  a\Psi^{(e)}_n(a f_{\theta,b})
 e^{\bm{j} 2\pi f_{\theta,b,2}\left[(\bm{r}_{-\theta} \bm{b})_2-c\cot(\theta)\right]}
A_2(-\frac{f_{\theta,b,2}}{\sin(\theta)})
\;d f_{\theta,b,2}\;d f_{\theta,b,1}\\
&=&\frac{1}{\left|\sin(\theta)\right|} \int_{-\infty}^{\infty}
\bm{j} \frac{f_{\theta,b,1}}{f_{\theta,b}}
e^{\bm{j}2\pi f_{\theta,b,1}\left[(\bm{r}_{-\theta} \bm{b})_1-c\right]} 
\tilde{A}_2(f_{\theta,b,1}^2;\bm{\xi}_2,c)\;df_{\theta,b,1},
\end{eqnarray*}
where the last equation defines $\tilde{A}_2(f_{\theta,b,1}^2;\bm{\xi}_2,c),$ as an {\em even} function of $f_{\theta,b,1},$ and
$\bm{f}_{\theta,b}= [f_{\theta,b,1}, f_{\theta,b,2} ] = \bm{r}_{-\theta} \bm{f}_b.$
Note that for fixed values of $a$ and $\theta=\theta_2$ we can find a value
of $\bm{b}$ such that
$(\bm{r}_{-\theta} \bm{b})_1=c,$ denoted $\bm{b}_2.$
When $\bm{\xi}=\bm{\xi}_2^{\prime}=\left[a,\theta_2,\bm{b}_2
\right]$ then $\sin{\left\{ (2 \pi f_{\theta,b,1})\left[(\bm{r}_{-\theta} \bm{b})_1-c\right]\right\} }$ vanishes identically for all $f_{\theta,b,1}.$
Hence $w_{n}^{(1)}(\bm{\xi}_2^{\prime};g_{s,2})=0,$
whilst
from equation
(\ref{rotenergy}) we find that the energy of the wavelet
transform with $\psi^{(q)}(\cdot)$ is conserved under rotations.
Thus
$\left|w^{(2)}(\bm{\xi};g)\right|$ is maximum at $\theta=\theta_2,$
and $\bm{b}=\bm{b}_2.$ 
With the correct choice of rotation $\theta$ we find
that we can retrieve the angle $\theta_2$ via maximising the energy of
$w_{n}^{(2)}(\bm{\xi};g_{s,2})$ and minimising 
the energy of $w_{n}^{(1)}(\bm{\xi};g_{s,2}).$  If the magnitudes of the
two wavelet transforms at $\bm{\xi}=\bm{\xi}_0$ are equal then $\theta=3\pi/4.$ Otherwise
we take a value of $\theta$ that maximises
\begin{eqnarray}
\nonumber
w_{n}^{(2)2}(\bm{\xi};g_{s,2})-w_{n}^{(1)2}(\bm{\xi};g_{s,2})
&=&\cos\left(2\theta\right)
\left(-w_{n}^{(1)2}(\bm{\xi}_0;g_{s,2})+w_{n}^{(2)2}(\bm{\xi}_0;g_{s,2})
\right)\\
&&-2 \sin\left(2\theta\right)w^{(1)}_n(\bm{\xi}_0;g_{s,2})
w_{n}^{(2)}(\bm{\xi}_0;g_{s,2}).
\end{eqnarray}
This has a stationary point at
\begin{equation}
\label{statmax}
\theta_{\max,n}=\frac{1}{2}\tan^{-1}\left(\frac{2w_{n}^{(1)}(\bm{\xi}_0;g_{s,2})
w_{n}^{(2)}(\bm{\xi}_0;g_{s,2})}
{w_{n}^{(1)2}(\bm{\xi}_0;g_{s,2})-w_{n}^{(2)2}(\bm{\xi}_0;g_{s,2})} \right),
\end{equation}
which corresponds to a maximum 
by choosing the appropriate solution.
Note that if $\left|w_{n}^{(2)}(\bm{\xi}_0;g_{s,2}\right|>\left|w_{n}^{(1)}(\bm{\xi}_0;g_{s,2})\right|$
we choose the solution $-\frac{\pi}{4}<\theta <\frac{\pi}{4}$ that corresponds
to a maximum
whilst if 
$\left|w_{n}^{(2)}(\bm{\xi}_0;g_{s,2})\right|<\left|w_{n}^{(1)}(\bm{\xi}_0;g_{s,2})\right|$
we take $\frac{\pi}{4}<\theta \le \frac{\pi}{2},$ or $-
\frac{\pi}{2}<\theta<-\frac{\pi}{4}.$
Thus at any fixed point $\bm{\xi}$ we can find the orientation that would
result from a line-discontinuity passing through $\bm{b},$ by utilising
the above equation, and this characterises local orientational structure.
Note also that at $\bm{\xi}=\bm{\xi}_2^{\prime}$ as this evaluates the wavelet
transform at values of $\bm{b}$ corresponding to the discontinuity, 
$\left|w_{n}^{(+)}(\bm{\xi}_2^{\prime};g_{s,2})\right|^2$ will be large.
A plot of 
parameter $\theta_{\max,n}$ should be combined with a plot of 
$\left|w_{n}^{(+)}(\bm{\xi}_2^{\prime};g_{s,2})\right|^2$
to verify that local rapid variation corresponding to an edge is present.

\section{The monogenic Wavelet Transform of AM/FM/OM Images}
Consider analysis of images of the form given by
equation (\ref{amfmom}). For more general classes of images, i.e. such as
images that are constrained to be positive we may add a constant term to
the model, but as noted by \cite{Gonnet}, the wavelet transform is a zero-mean
filter, and so this makes no difference to the subsequent analysis.
We then find that
\begin{eqnarray}
w_{n}^{(e)}(\bm{\xi};c_l)&=&\frac{1}{2}\int \int a_l(\bm{b})
\left[e^{2\bm{j}\pi 
\left(\phi_l(\bm{b}) +  \phi_l^{\prime}(\bm{b})
\bm{n}_l(\bm{b})\cdot(\bm{x}-\bm{b})\right)}+e^{-2\bm{j}\pi 
\left(\phi_l(\bm{b}) + \phi_l^{\prime}(\bm{b}) 
\bm{n}_l(\bm{b})\cdot(\bm{x}-\bm{b})\right)}\right]
\psi^{(e)*}_{\bm{\xi},n}(\bm{x})\;d^2\bm{x} \nonumber \\
&=&a_l(\bm{b})\cos[2\pi \phi_l(\bm{b})] 
a\Psi_{n}^{(e)}\left(a \phi_l^{\prime}(\bm{b}) \right), \nonumber \\
w_{n}^{(m)}(\bm{\xi}_0;c_l)&=&\frac{1}{2}\int 
\int a_l(\bm{b})
\left[e^{2\bm{j}\pi 
\left(\phi_l(\bm{b})+ \phi_l^{\prime}(\bm{b}) \bm{n}_l(\bm{b})
(\bm{x}-\bm{b})\right)}+
e^{-2\bm{j}\pi \left(\phi_l(\bm{b})+ 
\phi_{l}^{\prime}(\bm{b}) \bm{n}_l(\bm{b})(\bm{x}-
\bm{b})\right)}\right]\psi^{(m)*}_{\bm{\xi}_0,n}(\bm{x})\;d^2\bm{x} \nonumber \\
\label{quatterm}
&=&\frac{1}{2}a_l(\bm{b})\left[e^{2\bm{j}\pi 
\phi_l(\bm{b})}a\Psi^{(m)*}_{n}\left(a \phi_l^{\prime}(\bm{b})\bm{n}_l(\bm{b}) \right)
+e^{-2\bm{j}\pi 
\phi_l(\bm{b})}a\Psi^{(m)*}_{n}
\left(-a \phi_l^{\prime}(\bm{b})\bm{n}_l(\bm{b})\right)\right], \ m=1,2. \nonumber 
\end{eqnarray}
Hence it follows from equation (\ref{rotenergy}) that
$
w^{(+)}(\bm{\xi};c_l)=a_l(\bm{b})
a\Psi_{n}^{(e)}\left(a \phi_{l}^{\prime}(\bm{b})\right)e^{2\pi {\bm{e}}_{\nu_l}
\phi_l(\bm{b})},
$
with $\bm{e}_{\nu_l} = \bm{i} \cos\left(\eta_l(\bm{b})-\theta\right) + \bm{j} \sin
\left(\eta_l(\bm{b})-\theta\right)$.
This is the localised analogue of equation (\ref{monogenicloc}). For multi-component images, we may be able to separate the dominant
component out, similarly to ridge analysis based on complex wavelets \cite{Gonnet};
this requiring the assumption
$
a_l(\bm{b})
a\Psi_{n}^{(e)}\left(a \phi_{l}^{\prime}(\bm{b}) \right)>>
a_m(\bm{b})
a\Psi_{n}^{(e)}\left(a  \phi_{l}^{\prime}(\bm{b})  \right)\;\forall\;l\neq
m,
$
at all $\bm{\xi}$ considered. Furthermore,
the modulus of the wavelet transform is
\begin{equation}
\label{moddy2}
\left|w^{(+)}(\bm{\xi};c_l)\right|^2=a_l^2(\bm{b})
a^2\Psi_{n}^{(e)2}\left(a  \phi_{l}^{\prime}(\bm{b}) \right),
\end{equation}
and hence the wavelet transform of $c_l(\cdot)$ is locally maximal 
on the curve given by \\
$
{\mathcal{R}}(a,\theta,\bm{b})=\left\{(a,\bm{b}):\;a \phi_{l}^{\prime}(\bm{b})
=f_{\mathrm{max}}^{(n)}\right\},
$
where $f_{\mathrm{max}}^{(n)}$ is given by (\ref{fmaxxy}).
This defines the monogenic wavelet ridges \cite{Gonnet} of an AM/FM/OM images.
At any point
on this ridge, the local orientation may be computed. 
Ridge analysis is based on the fact that not all information of the redundant wavelet transform representation
needs to be considered to characterise the image:
only the ridge itself. As the ridge definition
does not depend on the angle $\theta,$ we need not carry out the transform
for all these values; a computational advantage to using directional wavelets. On the ridge we characterise the oscillatory
components locally as
\begin{eqnarray}
\nu_l(\bm{b})&=&\tan^{-1}\left(\frac{w^{(2)}(\bm{\xi};c_l)}
{w^{(1)}(\bm{\xi};c_l)}\right)=\eta_l(\bm{b})-\theta,\\
\phi_l(\bm{b})&=&\frac{1}{2\pi}\tan^{-1}\left(\frac{{\mathrm{sgn}}\left(w^{(2)}(\bm{\xi};c_l)
\right)\sqrt{w^{(1)2}(\bm{\xi};c_l)+
w^{(2)2}(\bm{\xi};c_l)}}
{w^{(e)}(\bm{\xi};c_l)}\right),\\
a^2_l(\bm{b})&=&  \frac{\left|w^{(+)}(\bm{\xi};c_l)\right|^2}{a^2\Psi^{(e)2}\left(
a \phi_{l}^{\prime}(\bm{b} )\right)}.
\end{eqnarray}
We have constrained $-\frac{\pi}{2}\le \nu_l \le
\frac{\pi}{2}$ and $-\frac{1}{2} \le \phi
\le \frac{1}{2}$ by the choice of sign for the $\tan^{-1}(\cdot).$

\section{Digital Implementation}
To preserve the exact monogenic
structure we implement the wavelet transform from the Fourier domain, calculating
the IDFT, thus making the algorithm of order $O(N_1 N_2 \log(N_1)\log(N_2)).$
We consider the maximum and minimum scales that can be resolved
- the range of the angle $\theta\in \left(0,2\pi\right)$ and $\bm{b}\in\left(0,N_1\Delta_1
\right)\times \left(0,N_2\Delta_2\right).$ As the real two-dimensional even
wavelet $\Psi_n^{(e)}(\cdot)$ is built from a real one-dimensional wavelet
corresponding to a band-pass filter, there exist frequencies $f_1^{(n)}$
and $f_2^{(n)}$ such that 
\begin{equation}
\Psi_n^{(e)}(f)\approx 0\;\forall f:\;\left|f\right| \notin\left(f_1^{(k)},f_2^{(n)}\right).
\end{equation}
Note that the DFT of observed image $g(\cdot,\cdot)$ is periodic by construction,
and that the standard assumption corresponds to $G(f_1,f_2)=0$ for all frequencies
not in the Nyquist band. Thus to perform the implementation we consider only scales $a$ such that
$
\Psi^{(e)}_n\left(a\sqrt{f_1^2+f_2^2}\right)\approx 0\;\forall\;\bm{f}:f_1\ge
1/(2\Delta_1),\;f_2\ge 1/(2\Delta_2).
$
This necessitates
$
a\ge a_{\min}=2 f_2^{(n)}\frac{\Delta_1 \Delta_2}{\sqrt{\Delta_1^2+\Delta_2^2}}.
$
As $a$ increases in magnitude, the wavelet becomes more peaked in the frequency
domain, and to ensure the wavelet covers at least $M$ frequency points we
constrain
$
a\le a_{\max}=\frac{1}{M}\min
\left(N_1\Delta_1,N_2\Delta_2\right)\left[f_2^{(n)}-f_1^{(n)}\right].
$

\section{Statistical Properties}
Consider estimation of features present in an image immersed in white noise
where the image is collected in a regular grid consisting of $x_1=s_1 \Delta_1,\;
s_1=0,\dots,N_1-1$
and $x_2=s_2 \Delta_2,\;s_2=0,\dots,N_2-1.$ We model the observed image 
$y(x_1 ,x_2 )$ as
$
y(x_1 ,x_2 )=g(x_1 ,x_2 )+\epsilon_{s_1,s_2}.
$
The noise $\left\{\epsilon_{s_1,s_2}\right\}$ is modelled as 
isotropically Gaussian and white. It is assumed that 
$E\left(\epsilon_{s_1,s_2}\right)=0,$
and $E\left(\epsilon_{s_1,s_2},\epsilon_{u_1,u_2}\right)=\delta_{s_1,u_1}
\delta_{s_2,u_2}\sigma^2_{\epsilon}.$
The wavelet transform of the noise will also be Gaussian, as it corresponds
to a sum of jointly normal variables. 
To give the distribution of the wavelet transform, we calculate its first
and second order structure at a fixed $\bm{\xi}.$ The wavelet transform is
a linear operation
and $w_{n}^{(\cdot)}(\bm{\xi};y)=
w_{n}^{(\cdot)}(\bm{\xi};g)+
w_{n}^{(\cdot)}(\bm{\xi};\epsilon).$
It follows that
$
E\left(w_{n}^{(\cdot)}(\bm{\xi};y)\right)= w_{n}^{(\cdot)}(\bm{\xi};g),
$
and we can determine the second order structure of the estimators from the
distribution of the noise.\\
We recast the full wavelet transform of the noise as a vector with
real valued entries, $
\bm{w}_{n}(\bm{\xi};\epsilon)
=\left[
w_{n}^{(e)}(\bm{\xi};\epsilon)\;\;
w_{n}^{(1)}(\bm{\xi};\epsilon)\;\;
w_{n}^{(2)}(\bm{\xi};\epsilon)
\right]^T,\;n=0,\dots,N-1.$
In Appendix B, we find with the additional assumption of
\begin{equation}
\min\left(\frac{1}{2\Delta_1},\frac{1}{2\Delta_2}\right)>\max_{n=0,\dots,N-1}
f_{\max}^{(n)},\end{equation}
that
\begin{equation}
\bm{w}_{n}(\bm{\xi};\epsilon)
\overset{d}{=}{\mathcal{N}}_3
\left(\bm{0}_3,\sigma^2_{\epsilon}\bm{V}\right),\;\bm{V}=
\begin{pmatrix}
1 & 0 & 0\\
0 & \frac{1}{2} & 0\\
0 & 0 & \frac{1}{2}
\end{pmatrix}. 
\end{equation}
\\
We perform estimation using the multiple orthogonal wavelets. Any
estimator of local signal properties needs to be smoothed, or averaged to
obtain
a low variance. The wavelet transform using any of the specified wavelet
functions averages the data across a window
in space and spatial frequency, where the width of the region depends on
the wavelet chosen, and in our case is characterised by the radial Morse
region ${\mathcal{D}},$ and the parameters $(\beta=l+\frac{1}{2},\gamma=m).$ Thomson \cite{Thomson}
suggested forming
estimates of local properties by averaging local energy estimates using several
orthogonal wavelets/functions. This usage explicitly reduces the
variability of the estimates with a clearly specified averaging region --
${\mathcal{D}}.$ Coherent behaviour over ${\mathcal{D}}$
is re-enforced across wavelet estimates, but the noisy uncorrelated
behaviour should cancel.
The bias inherent in the averaging is characterised
by the eigenvalues square
of the localisation operator.
In Appendix B we show that
$
E\left(w^{(l_1)}_{n_1}(\bm{\xi};\epsilon) 
w^{(l_2)}_{n_2}(\bm{\xi};\epsilon) \right)=
\sigma^2_{\epsilon} V^{l_1 l_2} \delta_{n_1,n_2}, 
$
and thus $\bm{w}_{n_1}(\bm{\xi};\epsilon)$ is uncorrelated with
$\bm{w}_{n_2}(\bm{\xi};\epsilon)$ unless $n_1=n_2,$ 
that combined with the assumption of Gaussian errors corresponds to independence.
We define averages of
the wavelet transform and the scalogram that will be used as a basis for
calculating estimators of other quantities as
\begin{equation}
\overline{w}^{(l)}(\bm{\xi};\cdot)=\frac{1}{N}\sum_{n=0}^{N-1}
w_n^{(l)}(\bm{\xi};\cdot)\;\;
\overline{S}^{(l)}(\bm{\xi};\cdot)=\frac{1}{N}\sum_{n=0}^{N-1}
S_n^{(l)}(\bm{\xi};\cdot),
\end{equation}
with $l=e,1,2,+,$ and finally as a measure of covariation we define
for $l_1,\;l_2=e,1,2,+,$ and images $g_1(\cdot),\;g_2(\cdot), 
\;\overline{C}^{(l_1,l_2)}(\bm{\xi};g_1(\cdot),g_2(\cdot))=
\frac{1}{N}\sum_{n=0}^{N-1}
w_n^{(l_1)}(\bm{\xi};g_1)w_n^{(l_2)*}(\bm{\xi};g_2).$
We define the estimators
$\widehat{w}^{(l)}(\bm{\xi};g) = \overline{w}^{(l)}(\bm{\xi};y)$ and
$\widehat{S}^{(l)}(\bm{\xi};g) = \overline{S}^{(l)}(\bm{\xi};y)$, for $l=e,1,2,+,$
as well as $\widehat{w^{(l_1)} w^{(l_2)}}(\bm{\xi};g_1)=
\overline{C}^{(l_1,l_2)}(\bm{\xi};g_1(\cdot),g_1(\cdot)).$
The Gaussian assumptions on $\epsilon$ then give 
$\overline{\bm{w}}(\bm{\xi};\epsilon)= 
[\overline{w}^{(e)}(\bm{\xi};\epsilon), \overline{w}^{(1)}(\bm{\xi};\epsilon), \overline{w}^{(2)}(\bm{\xi};\epsilon) ] 
 \overset{d}{=}
{\mathcal{N}}\left(\bm{0_3},\overline{\sigma}^{2}_{\epsilon} \bm{V}\right),$ where $\overline{\sigma}^{2}_{\epsilon}=\sigma^2_{\epsilon}/N.$
For most quantities we would intuitively expect to see a reduction of $1/N$
in their variances. 
When estimating the energy of the image at point $\bm{\xi}$
we consider 
\begin{eqnarray}
\widehat{S}^{(+)}(\bm{\xi};g)&=&\frac{1}{N} \sum_{n=0}^{N-1} 
S_n^{(+)}(\bm{\xi};y)\\
&=&
\overline{S}^{(+)}(\bm{\xi};g)+\overline{S}^{(+)}(\bm{\xi};\epsilon)+\frac{2}{N}
\sum_{n=0}^{N-1} 
\left[w_{n}^{(e)}(\bm{\xi};g)w_{n}^{(e)}(\bm{\xi};\epsilon)\right.\\
&&\left.+w_{n}^{(1)}(\bm{\xi};g)w_{n}^{(1)}(\bm{\xi};\epsilon)+
w_{n}^{(2)}(\bm{\xi};g)w_{n}^{(2)}(\bm{\xi};\epsilon)\right]
\end{eqnarray}
Up to
order $\overline{\sigma}^2_{\epsilon}$, we find, with the additional assumption
of the localised behaviour of $g()$ coherent across the $n$ wavelets,
\begin{equation}
\widehat{S}^{(+)}(\bm{\xi};g)\overset{d}{=}
{\mathcal{N}}\left(\overline{S}^{(+)}(\bm{\xi};g),4 \frac{\sigma^2_{\epsilon}}{N}
\left[\overline{S}^{(e)}(\bm{\xi};g)+\frac{1}{2} \left( \overline{S}^{(1)}(\bm{\xi};g) +  \overline{S}^{(2)}(\bm{\xi};g) \right)
\right]
\right).
\end{equation}
Hence the variance of the energy estimate decreases $O\left(\frac{1}{N}\right).$
\subsection{Distribution of Estimators}
We estimate the orientation of the line discontinuity in section VII by maximising
the difference
between the energy of the second and first components. Each wavelet indexed
by $n$ satisfies equation (\ref{statmax}) and thus 
writing the equations in terms of $\tan(\theta_{\max,n})=\theta_2,$ we may
sum over
the equations to find that
\[\theta_2=\frac{1}{2}\tan^{-1}\left(\frac{2 
\overline{C}^{(1,2)}(\bm{\xi};g_{s,2},g_{s,2}))
}
{\overline{S}^{(1)}(\bm{\xi}_0;g_{s,2})-\overline{S}^{(2)}(\bm{\xi}_0;g_{s,2})}\right).\]
We form estimate
\begin{equation}
\nonumber
\widehat{\theta}_{\rm{max}}(g_{s,2})=  
\frac{1}{2}\tan^{-1}\left(\frac{2 
\widehat{w^{(l_1)} w^{(l_2)}}(\bm{\xi};g_{s,2})
}
{\widehat{S}^{(1)}(\bm{\xi}_0;g_{s,2})-\widehat{S}^{(2)}(\bm{\xi}_0;g_{s,2})}\right).
\end{equation}
Let $w_{n}^{(l)}(\bm{\xi};\epsilon) = \sigma_{\epsilon}
w_{n,\epsilon}^{(l)}$, which entails that $\overline{w}^{(l)}(\bm{\xi};\epsilon)
= 
\sigma_{\epsilon} \overline{w}_{\epsilon}^{(l)}$
for $l=e,1,2,+$, and expand the above expression
$
\widehat{\theta}_{\rm{max}}(g_{s,2})= \theta_2 +
\sigma_{\epsilon} \delta \theta_2
+O(\sigma_{\epsilon}^{2}) 
.$
Note that $E\left(\delta \theta_2\right)=0$ making the estimator up to order
$\sigma^2_{\epsilon}$ unbiased, and the estimator has variance
\begin{equation}
{\mathrm{Var}}\left[ \widehat{\theta}_{\rm{max}}(g_{s,2}) \right]=
\frac{\sigma^{2}_{\epsilon}}{N}   \frac{1}{2}   \frac{
\overline{S}^{(1)}(g_{s,2}) + \overline{S}^{(2)}(g_{s,2})}
{\left[ \overline{S}^{(1)}(g_{s,2}) - \overline{S}^{(2)}(g_{s,2}) \right]^{2}
+ 4 \left[\overline{C}^{(1,2)}(g_{s,2},g_{s,2})\right]^2  } + O(\sigma^{3}_{\epsilon}).
\end{equation}
Thus, using multiple wavelets leads to a variance reduction.
\\
For AM/FM/OM signals we define the estimator for the orientation angle of the unit quaternion as
\begin{equation}
\widehat{\nu}_{l}(\bm{\xi};c_l)=\overline{\nu}_{l}(\bm{\xi};y_l)= \tan^{-1}\left[
\frac{\widehat{w}^{(2)}(\bm{\xi};c_l)}{\widehat{w}^{(1)}(\bm{\xi};c_l) }
\right],
\label{estnu}
\end{equation}
that up to order $\sigma^2_\epsilon$ is
\begin{equation}
\widehat{\nu}_{l}(\bm{\xi};c_l)
= 
\nu_{l}(\bm{\xi};c_l)+
\sigma_{\epsilon} \frac{\overline{w}^{(1)}(\bm{\xi};c_l) \overline{w}^{(2)}_{\epsilon}-
\overline{w}^{(2)}(\bm{\xi};c_l)\overline{w}^{(1)}_{\epsilon}}{\overline{S}^{(1)}(\bm{\xi};c_l)
+ 
\overline{S}^{(2)}(\bm{\xi};c_l)}+ O\left(\sigma_{\epsilon}^2\right)
\end{equation}
as
$\overline{w}^{(1)}(\bm{\xi}_0;c_l)=a_l(\bm{b}) \sin(2\pi \phi_l(\bm{b}))\cos(\nu_l)
\frac{a}{N}
\sum \Psi_{n_1}(a\left|\nabla \phi_l(\bm{b})\right| )$ and a similar results
holds for the second component
we have $\tan(\nu_l)=
\frac{\overline{w}^{(1)}(\bm{\xi};c_l)}{\overline{w}^{(2)}(\bm{\xi};c_l)}.$ As the wavelet transform of the noise has expectation zero
the estimator is unbiased and it has
variance
\begin{eqnarray}
{\rm Var}\left[\widehat{\nu}_{l}(\bm{\xi};c_l) \right] &=& \frac{ \sigma^{2}_{\epsilon}}{2N} \frac{\overline{w}^{(1)2}(\bm{\xi};c_l)+\overline{w}^{(2)2}(\bm{\xi};c_l) }{\left[ \overline{S}^{(1)}(\bm{\xi};c_l)
+\overline{S}^{(2)}(\bm{\xi};c_l)  \right]^{2}} + O(\sigma^{3}_{\epsilon})\\
&=&
\frac{ \sigma^{2}_{\epsilon}}{2N^2}\frac{(\sum \Psi_{n_1}(a\left|\nabla \phi_l(\bm{b})\right| ))^2}
{\sum \Psi_{n_1}^2(a\left|\nabla \phi_l(\bm{b})\right| ) }
\approx \frac{ \sigma^{2}_{\epsilon}}{2N}
.
\end{eqnarray}
Using multiple wavelets leads to variance reduction.
To estimate the phase we only use a single wavelet, the $n=0.$ Due to the orthogonality relations, the wavelet filters in the
Fourier domain cannot be strictly
positive for all frequencies, and thus for $n>0$ there are induced variations in the
phase estimate whenever
the wavelet filter changes sign. The Morse wavelet estimate of the phase
still profits from the wavelet's good
radial localisation.
\begin{eqnarray}
&& \widehat{\phi}_{l}(\bm{\xi};c_l) = \phi_{0,l}(\bm{\xi};c_l) + O\left(\sigma^2_{\epsilon}\right) \\
&& + 
\frac{\sigma_{\epsilon}}{2\pi}
\left[   \frac{-\left(S^{(1)}_0(\bm{\xi};c_l)+
S^{(2)}_0(\bm{\xi};c_l) \right) w^{(e)}_{\epsilon,0}              + w^{(e)}_{0}(\bm{\xi};c_l) \left[w^{(1)}_{\epsilon,0}
w^{(1)}_0(\bm{\xi};c_l)+w^{(2)}_{\epsilon,0} w^{(2)}_0(\bm{\xi};c_l)
\right]}{ S^{(+)}_{0}(\bm{\xi};c_l)     \sqrt{S^{(1)}_0(\bm{\xi};c_l)+
S^{(2)}_0(\bm{\xi};c_l)}  }\right] \nonumber.
\end{eqnarray}
As the expected value of the wavelet transform of noise is zero, the estimator
is thus unbiased, and
the variance of
the phase estimator is
\begin{equation}
{\mathrm{Var}}\left[\widehat{\phi}_{l}(\bm{\xi};c_l)
\right]
=\frac{\sigma^2_{\epsilon}}{S^{(+)}_0(\bm{\xi};c_l)} \frac{S^{(+)}_0(\bm{\xi};c_l)-\frac{1}{2} S^{(e)}_0(\bm{\xi};c_l)   }{(2\pi)^2 S^{(+)}_0(\bm{\xi};c_l)}
=\frac{\sigma^2_{\epsilon}}{(2\pi)^2 S^{(+)}_0(\bm{\xi};c_l)} \left[1-\frac{1}{2}\cos^{2}(2\pi \phi_l(\bm{b}))\right].
\end{equation}
When considering larger scales, the
wavelets are averaging across a lot of sample points, and the variance of
the phase estimate decreases. 
The amplitude is estimated as
$
\widehat{a}_{l}^{2}(\bm{b})= \frac{\widehat{S}^{(+)}(\bm{\xi};c_l)}{a^2 \Psi^{(e)2}(a
\widehat{\phi}_{l}^{\prime}(\bm{b})  )}. $

\section{Examples}
\begin{figure*}[t]
\centerline{
\includegraphics[height=1.75in,width=1.75in]{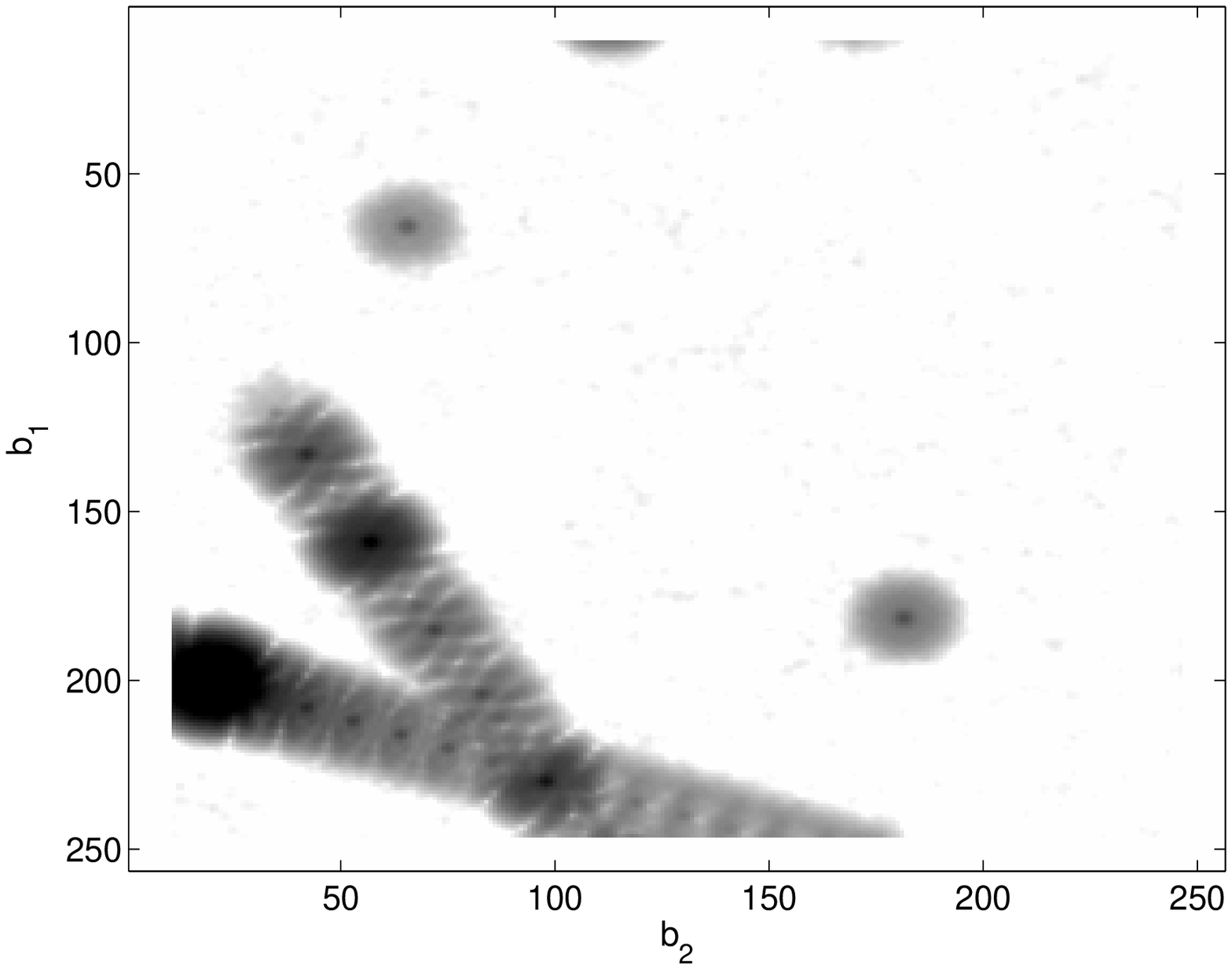}
\includegraphics[height=1.75in,width=1.75in]{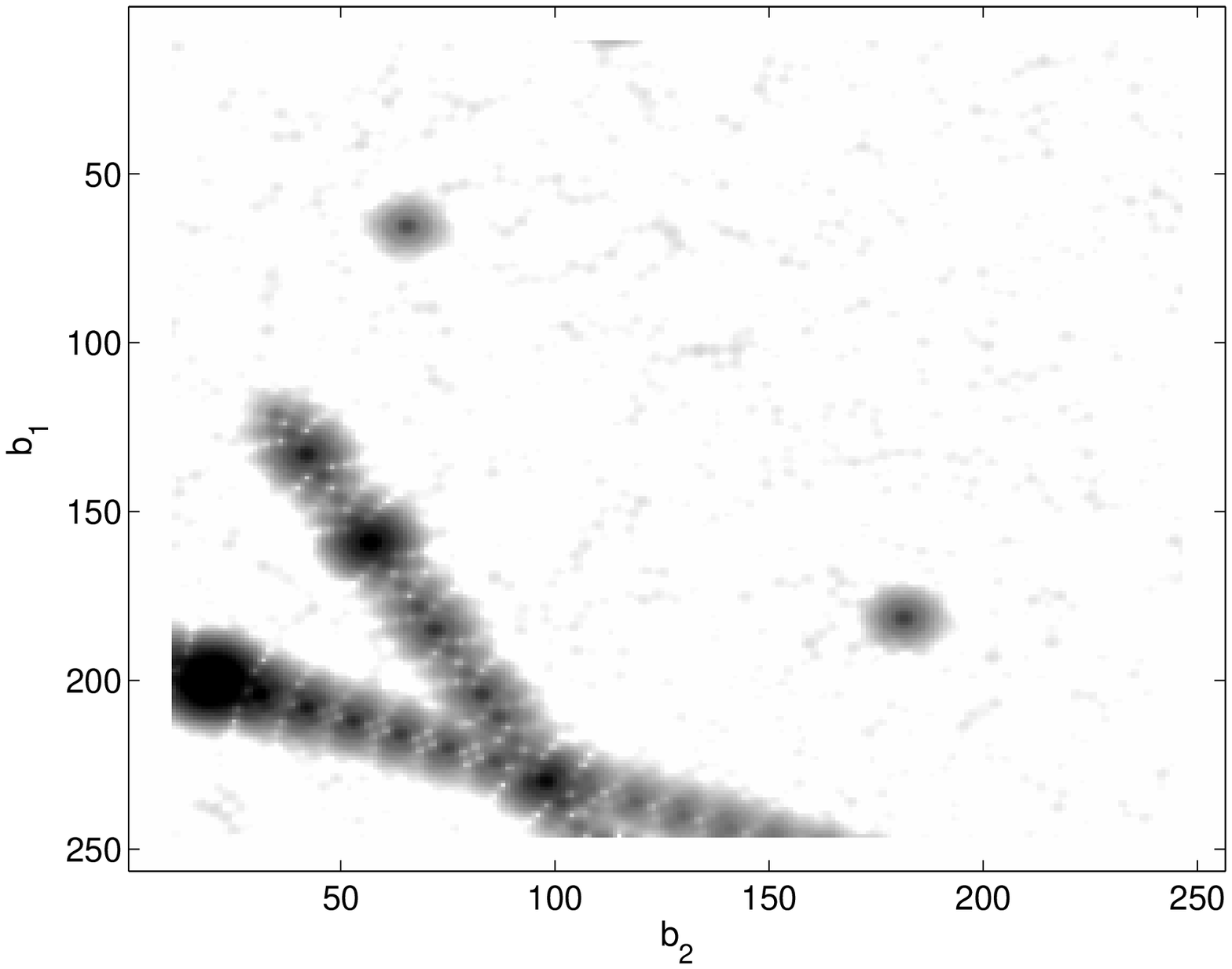}
\includegraphics[height=1.75in,width=1.75in]{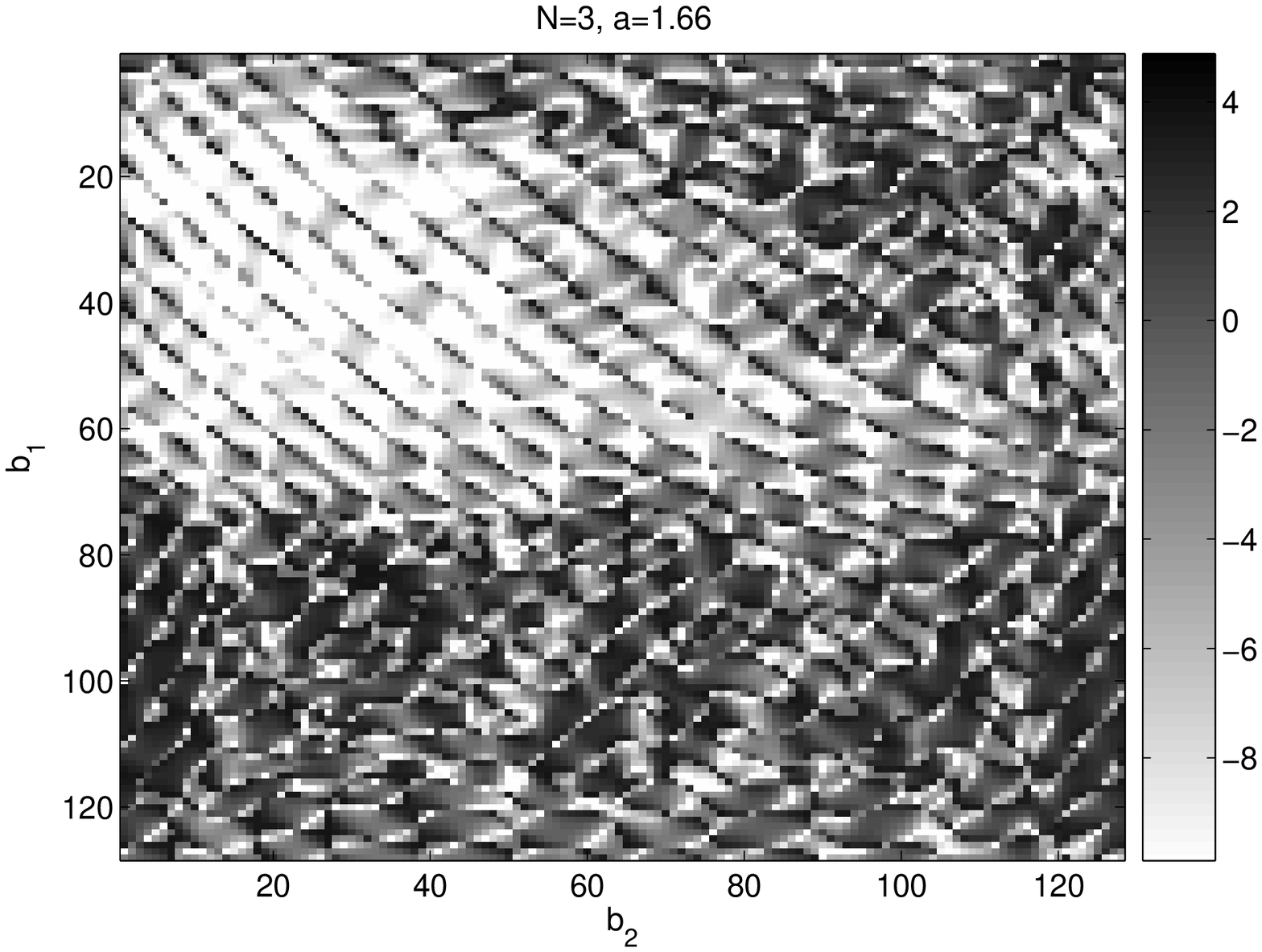}
\includegraphics[height=1.75in,width=1.75in]{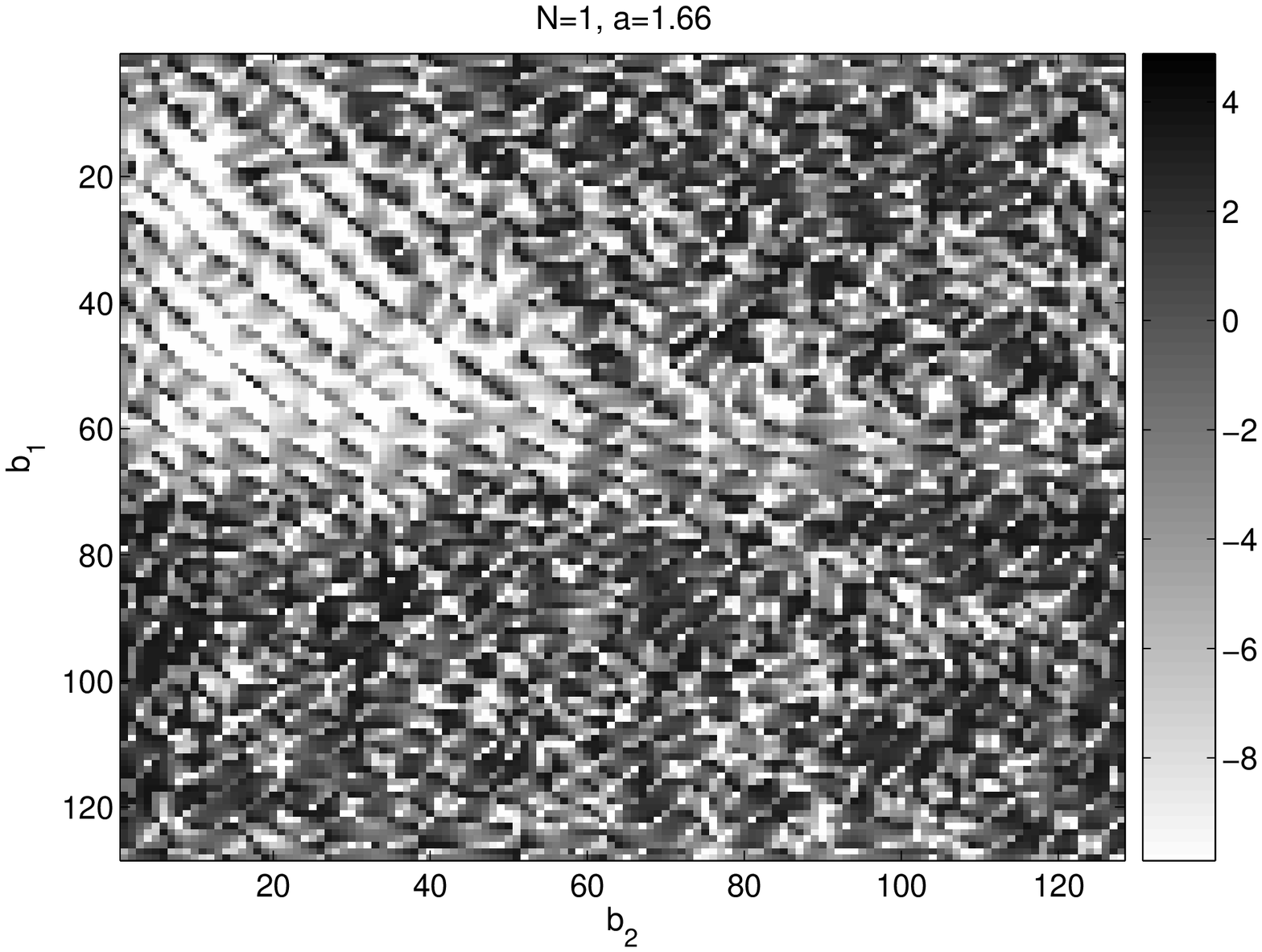}}
\caption{\label{fig:partition7} 
The local energy of signal 1 on a dB scale using
three wavelets (far left), or one wavelet (second from left), $a=1.4.$ 
The deviation of the estimated orientation from the true orientation 
at scale $a=1.66$ using dB scale, using three wavelets (second from right) and one
wavelet (far right). For the region where the signal has presence at those spatial
and scale points, the estimate using three wavelets is less noisy.
}
\end{figure*}

Consider
a collection of singularities observed in noise:
$
g_1(\bm{x})=\sum g_{2j}(\bm{x})+\sigma_1 \bm{\epsilon}_{
\bm{x}}$ where 
$g_{11}(\bm{x})=\frac{10}{\|\bm{x}-\left[\frac{1}{4}N_1+\frac{1}{2} \;\;
\frac{1}{4}N_1+\frac{1}{2}  \right]^T\|},$ 
$g_{12}(\bm{x})=\frac{15}{\|\bm{x}-\left[\frac{45}{64}N_1+\frac{1}{2}\;\;
\frac{45}{64}N_1+\frac{1}{2}\right]^T\|},$
$g_{13}(\bm{x})=
\left|x_1\cos(\pi/3)-x_2\sin(\pi/3)-\frac{15}{128}N_1+\frac{1}{2}\right|^{-1},$
$g_{14}(\bm{x})=\left|x_1\cos(\pi/9)-x_2\sin(\pi/9)-\frac{45}{64}N_1+\frac{1}{2} \right|^{-1},$
and we take $\sigma_1=0.2.$
The first two singularities of this signal are point
singularities and the latter two are line singularities. See Figure 
\ref{fig:partition7} for a plot of the scalogram
of the observed image at scale $a=1.4,$ corresponding to radial frequencies of $0.17.$
It can be clearly made out that the averaged estimate of the local energy
is a great deal more
robust to the noise. 
\\
Signal 2 is a multi component AM/FM/OM signal given by
$g_2(\bm{x})=g_{21}(\bm{x})+g_{22}(\bm{x})+\sigma^2_{\epsilon}\epsilon_{\bm{x}},$
where
$g_{21}(\bm{x})=1.2 I\left(x_1<N_1/2\right)\cos\left(
2\pi\times 0.087(\frac{t_2^2}{2N_1}+t_2)
\right),$
$t_l=x_1\cos(\eta_l)+x_2\sin(\eta_l),\;l=1,2,$
$\eta_1=-\frac{\pi}{4}+\frac{x_1+x_2-155}{10N_1},$
$g_{22}(\bm{x})=0.8
\cos\left(0.05\pi \left(\frac{t_2^2}{10N_1}+t_2\right)\right),$ 
$\eta_2=\frac{\pi}{5},$ and $N_1=N_2=128.$ 
We consider estimating its orientation at a
scale where the more rapid sinusoid is present near the left-hand side at
those frequencies, and find that our orientation estimate is substantively
less noise when using multiple wavelets, as is confirmed by  
Figure \ref{fig:partition7}. Clearly using the multiple wavelets is substantively
decreasing the variability of the estimator.

%

\section{Conclusions}
We show that the multiple monogenic Morse wavelets hold great potential for
digital image processing and analysis.
The monogenic Morse wavelets are the natural two-dimensional extension
of the analytic Morse wavelets, and are the eigenfunctions of a two-dimensional,
nonseparable, localisation operator. They form an orthogonal system, where the orthogonality
establishes the statistical properties of Gaussian noise. By averaging 
across  wavelets, estimators of local properties of the signal achieve reduced
variability. The monogenic
properties of the wavelets form a natural framework for determining local phase and orientation
properties. This framework explicitly parameterises the local orientation
of any variational structure
via a unit quaternion, and the localised analysis considers radial structures, thus using
the natural metric of Cartesian distances in the spatial domain. 

\section*{Acknowledgments}
SO \& GM would like to thank Dr Frederik Simons for valuable discussions
and SO would like to thank Professor Andrew Walden for introducing her to this
research
area. GM would like to acknowledge
the EPSRC (UK) for their financial support.

\appendix
\section*{A: Orthogonality Relations}
We know that if we choose $c=(2l+1)/m-1$ \cite[p.~2663]{Olhede2002}
then
\begin{eqnarray}
\int_{-\infty}^{\infty} \Psi^{(e)}_{n_1}\left(f\right)
\nonumber
\Psi^{(e)}_{n_2}\left(f\right)\;df&=&
2 A_{n_1;l,m} A_{n_2;l,m}
\int_0^{\infty} (2\pi \left|f\right|)^{2l} e^{-2(2 \pi \left|f\right|)^m}L_{n_1}^c\left(2
(2\pi \left|f\right|)^{m}\right)
L_{n_2}^c\left(2
(2\pi \left|f\right|)^{m}\right)\;df\\
\nonumber
&=& A_{n_1;l,m} A_{n_2;l,m}
\int_0^{\infty} (s/2)^{(2l+1)/m-1} e^{-s}L_{n_1}^c\left(s\right)
L_{n_2}^c\left(s\right)\;\frac{ds}{m}\\
&=&\delta_{n_1,n_2}.
\label{ortho1}
\end{eqnarray}
Consider the two-dimensional integral of two dimensional wavelets
\begin{eqnarray}
\nonumber 
\langle \psi^{(e)}_{n_1},\psi^{(e)}_{n_2}\rangle
&=&\frac{1}{\pi}
\int_{-\infty}^{\infty} \int_{-\infty}^{\infty}
A_{n_1;l,m} (2\pi f)^l e^{-(2 \pi f)^m}L_{n_1}^{c^{\prime}}\left(2 (2\pi
f)^{m}\right)
A_{n_2;l,m} (2\pi f)^l e^{-(2 \pi f)^m}L_{n_2}^{c^{\prime}}\left(2 (2\pi
f)^{m}\right)\;
df_1\;df_2\\
\nonumber
&=&2 \int_{0}^{\infty} 
A_{n_1;l,m} A_{n_2;l,m} (2\pi f)^{2l+1}e^{-2 (2 \pi f)^m}
L_{n_1}^{c^{\prime}}\left(2 (2\pi f)^{m}\right)
L_{n_2}^{c^{\prime}}\left(2 (2\pi f)^{m}\right)\;df\\
&=&A_{n_1;l,m} A_{n_2;l,m}
\nonumber
\int_0^{\infty}(s/2)^{(2l+2)/m-1}e^{-s}
L_{n_1}^{c^{\prime}}\left(s\right)
L_{n_2}^{c^{\prime}}\left(s\right)\;\frac{ds}{m}
\\
&=&\delta_{n_1,n_2}.
\label{ortho2}
\end{eqnarray}
Also note that
\begin{eqnarray*} 
\langle \psi^{(l)}_{n_1},\psi^{(l)}_{n_2}\rangle
&=&\frac{1}{\pi}
\int_{-\infty}^{\infty} \int_{-\infty}^{\infty}
A_{n_1;l,m} (2\pi f)^l e^{-(2 \pi f)^m}L_{n_1}^{c^{\prime}}\left(2 (2\pi
f)^{m}\right)\\
&&
A_{n_2;l,m} (2\pi f)^l e^{-(2 \pi f)^m}L_{n_2}^{c^{\prime}}\left(2 (2\pi
f)^{m}\right)\frac{f_l^2}{f_1^2+f_2^2}\;
df_1\;df_2,\;l=1,2.
\end{eqnarray*}
This implies that 
$
\langle \psi^{(1)}_{n_1},\psi^{(1)}_{n_2}\rangle+
\langle \psi^{(2)}_{n_1},\psi^{(2)}_{n_2}\rangle=
\langle \psi^{(e)}_{n_1},\psi^{(e)}_{n_2}\rangle=\delta_{n_1,n_2},
$
and as $\Psi^{(e)}_n(\bm{f})$ is radially symmetric we may deduce that
\begin{equation}
\langle \psi^{(1)}_{n_1},\psi^{(1)}_{n_2}\rangle=
\langle \psi^{(2)}_{n_1},\psi^{(2)}_{n_2}\rangle,
\end{equation}
and thus
\begin{equation}
\langle \psi^{(1)}_{n_1},\psi^{(1)}_{n_2}\rangle=\frac{1}{2}
\delta_{n_1,n_2}.
\label{ortho3}
\end{equation}
Finally, note that
\begin{eqnarray} 
\nonumber
\langle \psi^{(e)}_{n_1},\psi^{(1)}_{n_2}\rangle
&=&\frac{1}{\pi}
\int_{-\infty}^{\infty} \int_{-\infty}^{\infty}
A_{n_1;l,m} (2\pi f)^l e^{-(2 \pi f)^m}L_{n_1}^{c^{\prime}}\left(2 (2\pi
f)^{m}\right)\\
\nonumber
&&
A_{n_2;l,m} (2\pi f)^l e^{-(2 \pi f)^m}L_{n_2}^{c^{\prime}}\left(2 (2\pi
f)^{m}\right)(-i)\frac{f_1}{f_1^2+f_2^2}\;
df_1\;df_2\\
&=&0,
\label{ortho4}
\end{eqnarray} due to the integral of an odd function over a symmetric region
being zero. Similarly
\begin{equation} 
\langle \psi^{(e)}_{n_1},\psi^{(2)}_{n_2}\rangle
=0,\;
\langle \psi^{(1)}_{n_1},\psi^{(2)}_{n_2}\rangle
=0.
\label{ortho5}
\end{equation}

\section*{B: Calculation of Statistical Properties}
Define the discrete Fourier transform of the noise $\epsilon_{x_1,x_2}$,
\begin{equation}
\mathcal{E}(f_{1},f_{2})=\sum_{x_1=0}^{N_1-1}\sum_{x_2=0}^{N_2-1}
\epsilon_{x_1,x_2} e^{-2j\pi (f_1 x_1 + f_2 x_2)}.
\end{equation}
As the wavelet transform at any angle $\theta$ can be formed from
linear combinations of the wavelet transform at $\theta= 0$, in the way outlined
in section \ref{monor}, we need only calculate the properties at 
$\theta=0$. We have: 
\begin{eqnarray}
&&E\left(w_{n_1}^{(e)}(\bm{\xi};\epsilon)
\nonumber
w_{n_2}^{(e)*}(\bm{\xi};\epsilon)\right)\\
\nonumber
&=&\frac{a^2 \sigma^2_{\epsilon}}{N_1^2 N_2^2 \Delta_1^2 \Delta_2^2}
\sum_{l_1=-N_1^{\prime}}^{N_1^{\prime}-1}\sum_{l_2=-N_2^{\prime}}^{N_2^{\prime}-1}
\sum_{l_3=-N_1^{\prime}}^{N_1^{\prime}-1}\sum_{l_4=-N_2^{\prime}}^{N_2^{\prime}-1}
\delta_{l_1,l_3}\delta_{l_2,l_4} e^{2j\pi (b_1(l_1-l_3)+
b_2(l_2-l_4))}\\
\nonumber
&&\Psi_{n_1}\left(a\sqrt{\frac{l_1^2}{N_1^2 \Delta_1^2}+
\frac{l_2^2}{N_2^2 \Delta_2^2}}\right)
\Psi_{n_2}\left(a\sqrt{\frac{l_3^2}{N_1^2 \Delta_1^2}+
\frac{l_4^2}{N_2^2 \Delta_2^2}}\right)\\
&\approx& a^2 \sigma^2_1 \sigma^2_2\int_{-\frac{1}{2\Delta_1}}^{\frac{1}{2\Delta_1}}
\int_{-\frac{1}{2\Delta_2}}^{\frac{1}{2\Delta_2}}
\Psi_{n_1}\left(a\sqrt{f_1^2+f_2^2}\right)\Psi_{n_2}\left(a\sqrt{f_1^2+f_2^2}\right)\;
df_1\;df_2\\
&=&\sigma^{2}_{\epsilon} \delta_{n_1,n_2},
\end{eqnarray}
where the last line follows from Eq. (\ref{ortho2}). Similarly
\begin{eqnarray}
&&\cov\left(w_{n_1}^{(1)}(\bm{\xi};\epsilon), w_{n_2}^{(1)}(\bm{\xi};\epsilon) \right)\\
&\approx& a^2 \sigma_{\epsilon}^{2} \int_{-\frac{1}{2\Delta_1}}^{\frac{1}{2\Delta_1}} \int_{-\frac{1}{2\Delta_2}}^{\frac{1}{2\Delta_2}}
\Psi_{n_1}\left(a\sqrt{f_1^2+f_2^2}\right)\Psi_{n_2}\left(a\sqrt{f_1^2+f_2^2}\right)
\frac{f_1^2}{f_1^2+f_2^2}\;
df_1\;df_2\\
&=&\frac{1}{2}\sigma^2_{\epsilon}\delta_{n_1,n_2}.
\label{eqqy1}
\end{eqnarray}
following from Eq. (\ref{ortho3}), and likewise
\begin{eqnarray}
\cov\left(w_{n_1}^{(2)}(\bm{\xi};\epsilon),w_{n_2}^{(2)}(\bm{\xi};\epsilon) \right)
&=&\frac{1}{2}\sigma^2_{\epsilon}\delta_{n_1,n_2}. 
\label{eqqy12}
\end{eqnarray}

Furthermore,
\begin{eqnarray}
&&\cov\left(w_{n_1}^{(e)}(\bm{\xi};\epsilon),w_{n_2}^{(1)}(\bm{\xi};\epsilon) \right)\\
&=&E\left(w_{n_1}^{(e)}(a,0,b_1\Delta_1,b_2\Delta_2;\epsilon)
\nonumber
w_{n_2}^{(1)*}(a,0,b_1\Delta_1,b_2\Delta_2;\epsilon)\right)\\
\nonumber
&\approx& a^2 \sigma^2_{\epsilon} \int_{-\frac{1}{2\Delta_1}}^{\frac{1}{2\Delta_1}} \int_{-\frac{1}{2\Delta_2}}^{\frac{1}{2\Delta_2}}
\Psi_{n_1}\left(a\sqrt{f_1^2+f_2^2}\right)\Psi_{n_2}\left(a\sqrt{f_1^2+f_2^2}\right)
\frac{f_1}{\sqrt{f_1^2+f_2^2}}\;
df_1\;df_2\\
&=&0,
\label{eqqye1}
\end{eqnarray}
following from Eq. (\ref{ortho4}), likewise 
$\cov\left(w_{n_1}^{(e)}(\bm{\xi};\epsilon),w_{n_2}^{(1)}(\bm{\xi};\epsilon) \right)=0$ according to (\ref{ortho5}), 
 and finally
\begin{eqnarray}
&&\cov\left(w_{n_1}^{(1)}(\bm{\xi};\epsilon),w_{n_2}^{(2)}(\bm{\xi};\epsilon) \right)\\
&=&E\left(w_{n_1}^{(1)}(a,0,b_1\Delta_1,b_2\Delta_2;\epsilon)
\nonumber
w_{n_2}^{(2)*}(a,0,b_1\Delta_1,b_2\Delta_2;\epsilon)\right)\\
\nonumber
&\approx& a^2 \sigma^2_{\epsilon} \int_{-\frac{1}{2\Delta_1}}^{\frac{1}{2\Delta_1}} \int_{-\frac{1}{2\Delta_2}}^{\frac{1}{2\Delta_2}}
\Psi_{n_1}\left(a\sqrt{f_1^2+f_2^2}\right)\Psi_{n_2}\left(a\sqrt{f_1^2+f_2^2}\right)
\frac{f_1 f_2}{f_1^2+f_2^2}\;
df_1\;df_2\\
&=&0,
\label{eqqy13}
\end{eqnarray}
following from Eq. (\ref{ortho5}). This completes the covariance calculations for the distribution.

\end{document}